\documentclass[leqno]{amsart}

\usepackage{amsfonts}
\usepackage{amssymb, latexsym, amsmath, pb-diagram}
\usepackage{pstricks-add}
\usepackage{lamsarrow, pb-lams}
\usepackage{multicol}
\usepackage{bm}
\usepackage[mathscr]{eucal}
\usepackage{xy}
\usepackage{epic,eepic}
\xyoption{all}

\numberwithin{equation}{section}

\def\w{\widetilde}
\def\wh{\widehat}

\def\qb{\hfill $\Box$}
\def\nr{\refstepcounter{thm}\thethm}
\def\cf{{\it cf.}\ }

\newcommand{\PP}{\mathbb{P}}
\newcommand{\QQ}{\mathbb{Q}}

\newcommand{\kk}{\mathbf{k}}
\newcommand{\ba}{\mathbf{a}}
\newcommand{\xx}{\mathbf{x}}
\newcommand{\uu}{\mathbf{u}}
\newcommand{\vv}{\mathbf{v}}
\newcommand{\ZZ}{\mathcal{Z}}
\newcommand{\II}{\mathcal{I}}
\newcommand{\NN}{\mathbb{N}}

\begin{document}
\title[Simplicial complement and  Moment-angle Complexes]
{\bf The homology of simplicial complement and the cohomology of the moment-angle complexes}
\author[X. Wang \& Q. Zheng]{Xiangjun Wang and Qibing Zheng}
\thanks{The authors were supported by NSFC grant No. 10771105 and No. 11071125}
\address{School of Mathematical Sciences and LPMC, Nankai University,
Tianjin, 300071, P.R.China.}
\email{xjwang@nankai.edu.cn}
\address{School of Mathematical Science and LPMC, Nankai University,
Tianjin, 300071, P.R.China}
\email{zhengqb@nankai.edu.cn}
\maketitle

\begin{abstract}
A simplicial complement $\PP=\{\sigma_1, \sigma_2, \cdots , \sigma_s\}$ is a sequence
of subsets of $[m]$ and the simplicial complement $\PP$ corresponds to
an unique simplicial complex $K_\PP$ with vertices in $[m]$. In this paper, we defined
the homology of a simplicial complement $H_{i,\sigma}(\Lambda^{*,*}[\PP], d)$ over a principle ideal domain
$\kk$ and proved that $H_{*,*}(\Lambda[\PP], d)$ is isomorphic to
the $Tor$ of the corresponding face ring $\kk(K_\PP)$ by the Taylor resolution.  As applications, we give methods
to compute the ring structure of $Tor_{*,*}^{\kk[\xx]}(\kk(K_\PP), \kk)$,
$\mbox{link}_{K_\PP}\sigma$, $\mbox{star}_{K_\PP}\sigma$ and the cohomology modules of the generalized moment-angle complexes.
\end{abstract}

\section{Introduction and statement of results}

   The moment-angle complexes have been studied by topologists for many years
(\cf \cite{Porter} \cite{Lopez}). In 1990's
Davis and Januszkiewicz \cite{DJ} introduced quasi-toric manifolds which were being studied
intensively by algebraic geometers. They observed
that every quasi-toric manifold is the quotient of a moment-angle complex by the free action
of a real torus, here the moment-angle complex is denoted by $\mathcal{Z}_K$ corresponding to
an abstract simplicial complex $K$.
The topology of $\ZZ_K$ is complicated and getting more attentions by topologists lately
(\cf \cite{GM} \cite{Hochster} \cite{Baskakov} \cite{Panov} \cite{Franz}).
Recently a lot of work has been done on generalizing the moment-angle complex
$\ZZ_K=\ZZ_K(\underline{D^2}, \underline{S^1})$ to pairs of spaces $(\underline{X}, \underline{A})$
(\cf  \cite{BB}, \cite{BBCG08}, \cite{GTh}, \cite{LP}).
In this paper we study the cohomology of the generalized moment-angle complexes
$\ZZ_K(\underline{X}, \underline{A})$ corresponding to the pairs of spaces
$(\underline{X}, \underline{A})$ with inclusions $A_i\hookrightarrow X_i$ being
homotopic to constant for all $i$.

    Classically the homological algebra aspect of the Stanley-Reisner {\it face ring}
plays an important role in the cohomology of $\ZZ_K$ (\cf \cite{BP}, \cite{Stanley}).
Let $K$ be an abstract simplicial complex with $m$ vertices.
Choose a ground ring $\kk$ with unit (we are mostly interested in a
principal ideal domain).  Let $\kk[\xx]$ be the $\NN^m$ graded polynomial algebra
over $\kk$ on $m$ indeterminates $\xx=\{x_1, x_2, \cdots, x_m\}$. The Stanley-Reisner {\it face ring}
is the quotient ring $\kk(K)=\kk[\xx]/\II_K$, where $\II_K$ is the Stanley-Reisner ideal
generated by the monomials $\xx_\sigma$
corresponding to non-faces $\sigma\not\in K$.
\[
\II_K=\langle \xx_\sigma=x_{i_1}x_{i_2}\cdots x_{i_n}
    |\sigma=\{i_1, i_2, \cdots, i_n\}\not\in K\rangle.
\]
It is well known that the cohomology of $\ZZ_K$  is isomorphic to the $Tor$  of
$\kk(K)$ over $\kk[\xx]$, that is
\begin{align*}
H^{q}(\ZZ_K, \kk)\cong & \bigoplus_{2j-i=q} Tor^{\kk[\xx]}_{i, j}(\kk(K), \kk) & & \mbox{and}\\
 Tor^{\kk[\xx]}_{i, j}(\kk(K), \kk) \cong & Tor^{\kk[\xx]}_{i, j}(\kk, \kk(K)).
\end{align*}

    There are several methods to compute the $Tor$ of $\kk(K)$. For example, one can use the
Koszul resolution $\Lambda[u_1, u_2, \cdots, u_n]\otimes \kk[\xx]$ of $\kk$, then by applying
the functor $\otimes_{\kk[\xx]}\kk(K)$, the homology of $\Lambda[u_1, u_2, \cdots, u_n]\otimes \kk(K)$
is $Tor^{\kk[\xx]}_{i,2j}(\kk, \kk(K))$. On the other hand, one can construct a free resolution
$(R^*, d)$ of $\kk(K)$, then by applying the functor $\otimes_{\kk[\xx]}\kk$, the homology of $R^*\otimes_{\kk[\xx]}\kk$
is $Tor^{\kk[\xx]}_{i,j}(\kk(K), \kk)$. In this paper, we will follow the second method to
give a free resolution of $\kk(K)$ and  compute the homology of the resulting chain complex
$R^*\otimes_{\kk[\xx]}\kk$.

   Let $\PP=\{\sigma_1, \sigma_2, \cdots , \sigma_s\}$  be a sequence of subsets of
$[m]=\{1,2,\cdots,m\}$, which is called a {\it simplicial complement} in this paper. It is easy to see that
$$\II_\PP=\langle \xx_{\sigma_1}, \xx_{\sigma_2}, \cdots, \xx_{\sigma_s}\rangle$$
is an ideal generated by square-free monomials. Thus there is an unique simplicial complex $K_\PP$
such that $\II_{K_\PP}=\II_\PP$.

    Let $\Lambda[\PP]$ be the exterior algebra over $\kk$
generated by $\PP$. Given a generator $\uu=\sigma_{k_1}\sigma_{k_2}\cdots \sigma_{k_q}$ in $\Lambda[\PP]$,
let $S_\uu = \sigma_{k_1}\cup \sigma_{k_2}\cup \cdots \cup \sigma_{k_q}$ be the total subset
of $\uu$. The bi-degree of $\uu$ is defined as
 \[bideg(\uu)=(q, S_\uu),\]
and the monomial
$\xx_{S_\uu}$ determines an unique $\NN^m$ vector $\ba_\uu$ such that $\xx_{S_\uu}=\xx^{\ba_\uu}$.

Let
\[
\partial_i(\uu)
  =\sigma_{k_1}\cdots\widehat{\sigma}_{k_i}\cdots\sigma_{k_q}=
  \sigma_{k_1}\cdots\sigma_{k_{i-1}}\sigma_{k_{i+1}}\cdots\sigma_{k_q}
\]
and for a generator $\uu\otimes\xx^\ba\in \Lambda[\PP]\otimes_\kk\kk[\xx]$ define
\[
\widetilde{d}(\uu\otimes\xx^\ba)=\sum_i(-1)^i\partial_i(\uu)\otimes\xx_{(S_\uu\setminus S_{\partial_i(\uu)})}\xx^\ba,
\]
where $\displaystyle{\xx_{(S_\uu\setminus S_{\partial_i(\uu)})}=\frac{\xx_{S_\uu}}{\xx_{S_{\partial_i(\uu)}}}}$
denote the monomial of $\kk[\xx]$ corresponding to
the subset $S_\uu\setminus S_{\partial_i(\uu)}$.

   It is known that $(\Lambda[\PP]\otimes_\kk\kk[\xx], \w{d})$ is a free
$\kk[\xx]$-resolution of $\kk[\xx]/\II_\PP$  called the Taylor resolution \cite{Mil}.
Then applying the functor $\otimes_{\kk[\xx]}\kk$, we get
$$(\Lambda[\PP]\otimes_\kk\kk[\xx]\otimes_{\kk[\xx]}\kk, \w{d}\otimes 1)=(\Lambda^{*,*}[\PP],d)$$
and

  {\bf Theorem 2.7} {\it Let $\PP=\{\sigma_1, \sigma_2, \cdots, \sigma_s\}$ be a simplicial complement and $K_\PP$
be the corresponding simplecial complex such that $\II_{K_\PP}=\II_\PP$.
Let $(\Lambda^{*,*}[\PP], d)$ be the chain complex induced from the bi-graded exterior algebra on
$\PP$. The differential
$d: \Lambda^{q,*}[\PP]\rightarrow \Lambda^{q-1,*}[\PP]$ is given by
\begin{align*}
d(\uu)= & \sum_i(-1)^i\partial_i(\uu)\cdot\delta_{\partial_i(\uu)}
\end{align*}
where $\delta_{\partial_i(\uu)}=1$ if the total subset $S_\uu=S_{\partial_i(\uu)}$
and $\delta_{\partial_i(\uu)}=0$ if $S_\uu\not=S_{\partial_i(\uu)}$.  Then the homology of
$(\Lambda^{*,*}[\PP],\ d)$ is
\[
H_{*,\sigma}(\Lambda^{*,*}[\PP],\ d)=Tor^{\kk[\xx]}_{*,\sigma}(\kk[\xx]/\II_\PP, \kk)
  =Tor^{\kk[\xx]}_{*,\sigma}(\kk(K_\PP), \kk),
\]
which is called the homology of the simplicial complement $\PP$.}

    Fix a $\sigma\in 2^{[m]}$ and let $\Lambda^{*,\sigma}[\PP]$ be the submodule generated by
$\uu$ with total subset $S_\uu=\sigma$. From the formula for $d(\uu)$, we see that
$(\Lambda^{*,\sigma}[\PP], d)$ is a  sub-complex of
$(\Lambda^{*,*}[\PP], d)$. Then from a theorem of Baskakov, we see that
\[
H_{q,\sigma}(\Lambda^{*,\sigma}[\PP], d)=Tor^{\kk[\xx]}_{q,\sigma}(\kk(K_\PP), \kk)
  =Tor^{\kk[\xx]}_{q,\sigma}(\kk, \kk(K_\PP))=\w{H}^{|\sigma|-q-1}(K_\PP\cap\sigma, \kk),
\]
where $K_\PP\cap\sigma$ is the {\it full sub-complex} of $K_\PP$ consisting
of all simplices of $K_\PP$ which have all of their vertices in $\sigma$, that is
$K_\PP\cap\sigma=\{\tau\cap\sigma|\tau\in K_\PP\}$.

    Notice that both $\kk(K_\PP)$ and $\kk$ are $\kk[\xx]$ algebras, $Tor_{*,*}^{\kk[\xx]}(\kk(K_\PP), \kk)$
has natural algebraic structure and the following isomorphism of algebras holds.
\[
H^*(\ZZ_{K_\PP}, \kk)\cong Tor_{*,*}^{\kk[\xx]}(\kk(K_\PP), \kk).
\]
In section 3, we proved that

{\bf Theorem 3.6} {\it The algebraic structure in $H_{*,*}(\Lambda^{*,*}[\PP], d)\cong Tor_{*,*}^{\kk[\xx]}(\kk(K_\PP), \kk)$
is given by
\[
[c]\times[c']=[c\cdot c']\cdot \delta_{\sigma, \sigma'}
\]
where $[c]\in H_{q,\sigma}(\Lambda^{*,\sigma}[\PP], d)$, $[c']\in H_{j,\sigma'}(\Lambda^{*,\sigma'}[\PP], d)$,
the cycle $c\cdot c'$ is the product of $c$ and $c'$ in the exterior algebra $\Lambda^{*,*}[\PP]$ and
$\delta_{\sigma, \sigma'}=1$ if $\sigma\cap\sigma'=\phi$, $\delta_{\sigma, \sigma'}=0$ if $\sigma\cap\sigma'\not=\phi$.
}

    Fix a subset $\omega\subset [m]$, there is the {\it $\omega$-compression $E_\omega$} defined on the
simplicial complement $\PP$ given by
$$E_\omega\PP=\{\sigma_1\setminus\omega, \sigma_2\setminus\omega, \cdots, \sigma_s\setminus\omega\}.$$
For the subset $\omega$, the {\it link} and {\it star} of the simplicial complex $K_\PP$ corresponding to $\omega$
are the sub-complexes
\begin{align*}
\mbox{star}_{K_\PP}\omega= & \{\tau\in K_\PP|\tau\cup\omega\in K_\PP\}\\
\mbox{link}_{K_\PP}\omega= & \{\tau\in K_\PP|\tau\cup\omega\in K_\PP;\ \tau\cap\omega=\phi\}.
\end{align*}
In section 4 we proved that:

  {\bf Theorem 4.7} {\it Let $\PP$ be a simplicial complement and $K_\PP$ be the simplicial complex corresponding
to $\PP$. Then:}
\begin{enumerate}
\item {\it The $star$,} $\mbox{star}_{K_\PP}\omega=K_{E_\omega\PP}$ {\it is the simplicial
complex corresponding to the simplicial complement $E_\omega\PP$. Thus}
\[
Tor^{\kk[\xx]}_{q,\sigma}(\kk(\mbox{star}_{K_\PP}\omega), \kk)=
    H_{q,\sigma}(\Lambda^{*,\sigma}[E_\omega\PP], d),
\]
\item {\it The $link$,} $\mbox{link}_{K_\PP}\omega=K_{E_\omega\PP}\cap ([m]\setminus\omega)$ {and}
\[
H_{q,[m]\setminus \omega}(\Lambda^{*,[m]\setminus \omega}[E_\omega\PP], d)=
\w{H}^{m-|\omega|-q-1}(\mbox{link}_{K_\PP}\omega, \kk).
\]
\end{enumerate}

  Let $(\underline{X}, \underline{A})=\{(X_i, A_i, x_i)|i\in[m]\}$ denote
a set of triples of $CW$-complexes with base points $x_i\in A_i$ and $K_\PP$ be an abstract simplicial
complex. The generalized moment-angle complex determined by $(\underline{X}, \underline{A})$ and
$K_\PP$ denoted by $\ZZ_{K_\PP}(\underline{X}, \underline{A})$ is defined
to be the $colimit$
$$\ZZ_{K_\PP}(\underline{X}, \underline{A})=
\bigcup_{\sigma\in K_\PP} D(\omega),
$$
where
\begin{align*}
D(\omega)= & Y_1\times Y_2 \times \cdots \times Y_m  &\mbox{and} & &
Y_i=& \left\{\begin{array}{lll}
          X_i & \mbox{if $i\in\omega$}\\
          A_i & \mbox{if $i\not\in\omega$.}
            \end{array}\right.
\end{align*}

    Based on Theorem 4.7 and the decomposition of $\Sigma\ZZ_{K_\PP}(\underline{X}, \underline{A})$ given by
Bahri, Bendersky, Cohen and Gitler (\cf \cite{BB,BBCG08}), we proved in section 5

 {\bf Theorem 5.7} {\it Let $K_\PP$ be an abstract simplicial complex corresponding to a simplicial complement $\PP$
and let $(\underline{X}, \underline{A})=\{(X_i, A_i, x_i)|i\in [m]\}$ denote m choices of
connected, pointed pairs of $CW$-complexes, with the inclusion $A_i\hookrightarrow X_i$ homotopic to constant for
all $i$. If $\w{H}^*(X_i, \kk)$ and $\w{H}^*(A_i, \kk)$ are free $\kk$-modules for all $i$,
then the cohomology of $\ZZ_{K_\PP}(\underline{X}, \underline{A})$ is isomorphic to
\[
H^*(\ZZ_{K_\PP}(\underline{X}, \underline{A}), \kk)=
\bigoplus_{\omega\in K_\PP}\left(\bigoplus_{\tau\subset [m]\setminus\omega}
  H_{*,\tau}(\Lambda^{*,\tau}[E_\omega\PP], d) \otimes
  \left(\bigotimes_{i\in\omega}\w{H}^*(X_i, \kk)\right)\otimes
  \left(\bigotimes_{j\in\tau}\w{H}^*(A_j, \kk)\right)\right)
\]
as $\kk$-modules, where
\[
H_{*,\tau}(\Lambda^{*,\tau}[E_\omega\PP],d)\cong Tor_{*,\tau}^{\kk[\xx]}(\kk(\mbox{star}_{K_\PP}\omega), \kk).
\]
}

    As applications, we also consider the cohomology of $\ZZ_{K_\PP}(\underline{X}, \underline{A})$
for some special triples of $CW$-complexes $(\underline{X}, \underline{A})$ including
all the $X_i$ are contractible; all the $A_i$ are contractible and
$(\underline{X}, \underline{A})=(\underline{S^2}, \underline{S^1})$.

{\it Acknowledgements } The authors are grateful to Z.~L\"u for his introduction of this topic. The authors
are indebted to V.~M.~Buchstaber for his comments and suggestions, especially for the name of definition
{\it simplicial complement} which is originally called {\it partition.} The authors also thank S. Gitler for his
helpful suggestions.

\section{The homology of simplicial complement and the $Tor$ of face ring}

    Let $[m]=\{1,2,\cdots, m\}$ and $2^{[m]}$ denote the power set of $[m]$.
Choose a ground ring $\kk$ with unit (We usually suppose that $\kk$ is a principle ideal domain).
Let $\kk[\xx]$ be the $\NN^m$ graded polynomial algebra over $\kk$ on $m$ indeterminates
$\xx=\{x_1, x_2, \cdots, x_m\}$. The monomial of $\kk[\xx]$ is expressed as
\[
\xx^\ba=x_1^{a_1}x_2^{a_2}\cdots x_m^{a_m}
\]
for a unique vector $\ba=(a_1, a_2, \cdots, a_m)\in \NN^m$.
For a subset $\sigma=\{i_1, i_2, \cdots, i_n\}$ of $[m]$,
let $\xx_\sigma=x_{i_1}x_{i_2}\cdots x_{i_n}$ be the corresponding monomial in $\kk[\xx]$,
while $\xx_\phi=1$. Then the monomial
$\xx_\sigma$ can be expressed as $\xx^{\ba_\sigma}=\xx_\sigma$ for a unique
vector $\ba_\sigma\in\NN^m$.

  {\bf Definition \nr} {\it A sequence $\PP=\{\sigma_1, \sigma_2, \cdots , \sigma_s\}$ of
subsets of $[m]$ is called a simplicial complement. Given a simplicial complement
$\PP=\{\sigma_1, \sigma_2, \cdots, \sigma_s\}$,
let $\II_\PP$ denote the ideal of $\kk[\xx]$ generated by $\xx_{\sigma_i}$'s
\[
\II_\PP=\langle \xx_{\sigma_1}, \xx_{\sigma_2}, \cdots, \xx_{\sigma_s}\rangle.
\]
Two simplicial complements $\PP$ and $\QQ$  are called equivalent if they generate the same ideal of $\kk[\xx]$,
that is $\II_\PP=\II_\QQ$.}

    Given an abstract simplicial complex $K$, let $\PP_K$ denote the simplicial complement
consists of all the non-faces $\sigma\not\in K$ or equivalently consists of all the
{\it missing faces} in the sense that it is not a simplex of $K$ but all of its proper subsets
are simplices of $K$.  One can easily see from its definition
that $\II_{\PP_K}=\II_K$.

    Furthermore, given a simplicial complement $\PP$, noticed that $\II_\PP$ is
an ideal generated by square-free monomials,
there is a unique simplicial complex $K_\PP$ such that $\II_\PP=\II_{K_\PP}$.

    {\bf Definition \nr} {\it Define $K_\PP$  to be the simplicial complex corresponding to
$\PP$ such that the non-faces $2^{[m]}\setminus K_\PP$ is the
full subset of $2^{[m]}$ consisting of all subsets $\sigma\in 2^{[m]}$ which contain a
subset $\sigma_i$ in $\PP$, that is
\[
2^{[m]}\setminus K_\PP=\{\sigma\in 2^{[m]} \mid\ \mbox{there is a $\sigma_i\in\PP$
 such that $\sigma_i\subset\sigma$}\}.
\]}

     Notice that,
$\xx_\sigma=\xx_{\sigma_i}\cdot\xx_{\sigma\setminus\sigma_i}$ if $\sigma_i\subset\sigma$,
we have $\II_\PP=\II_{K_\PP}$.
The Stanley-Reisner face ring is
\[
\kk(K_\PP)=\kk[\xx]/\II_{K_\PP}=\kk[\xx]/\II_\PP.
\]

    {\bf Proposition \nr} {\it If $\PP$ and $\QQ$ are equivalent simplicial complements,
then
\[
\kk(K_\PP)=\kk[\xx]/\II_\PP=\kk[\xx]/\II_\QQ=\kk(K_\QQ)
\]
and the corresponding simplicial complexes $K_\PP=K_\QQ$. \qb}

    {\bf Remark:}  Notice that if the empty set $\phi$ is an element of the simplicial complement $\PP$, then
$\II_\PP=\langle \xx_\phi=1\rangle=\kk[\xx]$, the simplicial complex $K_\PP=\varnothing$, which is different from the
empty simplicial complex $K=\{\phi\}$.

{\bf Definition \nr} {\it   Given a simplicial complement
$\PP=\{\sigma_1, \sigma_2, \cdots , \sigma_s\}$,
let $\Lambda[\PP]$  be the exterior algebra over $\kk$ on the generators $\sigma_1, \sigma_2, \cdots , \sigma_s$.
Define an $\NN\times 2^{[m]}$-graded $\kk$-module structure on $\Lambda [\PP]$ as follows.
For a generator $\uu=\sigma_{k_1}\sigma_{k_2}\cdots\sigma_{k_q}$ of $\Lambda^{*,*}[\PP]$,
let
\begin{align*}
S_\uu = & \sigma_{k_1}\cup \sigma_{k_2}\cup \cdots \cup \sigma_{k_q}
& \mbox{and}& & S_1= & \phi
\end{align*}
denote the union of the subsets $\sigma_{k_j}$'s, which is called the {\it total subset} of $\uu$.
Define the degree of $\uu=\sigma_{k_1}\sigma_{k_2}\cdots\sigma_{k_q}$ as
\[
bideg(\uu)=(q, S_\uu)
\]
so that $\Lambda[\PP]$ becomes an $\NN\times 2^{[m]}$ graded $\kk-$module
and $\displaystyle{\Lambda^{*,*}[\PP]=\bigoplus_{q,\sigma}\Lambda^{q,\sigma}[\PP]}$. }

    {\bf Remark:}  $\xx_{S_\uu}$ is the least common multiple of the monomials
$\xx_{\sigma_{k_1}}, \xx_{\sigma_{k_2}}, \cdots, \xx_{\sigma_{k_q}}$ and
$\Lambda^{*,*}[\PP]$ corresponding to the $\xx_{\sigma_1},\ \xx_{\sigma_2},
\ \cdots,\ \xx_{\sigma_s}$ labeled standard $s-1$ simplex $\Delta^{s-1}$ (\cf \cite{Mil} Chapter 4).

    Let $R^{*,*}=\Lambda^{*, *}[\PP]\otimes_\kk \kk[\xx]$ denote the free $\kk[\xx]$-module
generated by $\Lambda^{*,*}[\PP]$. For a generator $\uu=\sigma_{k_1}\sigma_{k_2}\cdots\sigma_{k_q}$
of $\Lambda^{q,*}[\PP]$, the monomial in $\kk[\xx]$ corresponding to the total subset
$S_\uu$ is expressed as $\xx_{S_\uu}=\xx^{\ba_{\uu}}$ for an unique $\NN^m$ vector $\ba_\uu\in \NN^m$.
We define the $\NN\times\NN^m$ degree of  $\uu$ as
$$
bideg(\uu)=(q, \ba_{\uu}).
$$
Then $R^{*,*}$ becomes an $\mathbb{N}\times\mathbb{N}^m$ graded module over $\kk$.

   Given a generator $\uu=\sigma_{k_1}\sigma_{k_2}\cdots\sigma_{k_q}$
of $\Lambda^{q, *}[\PP]$, denote that
\begin{equation*}
\partial_i(\uu)=\partial_i(\sigma_{k_1}\sigma_{k_2}\cdots\sigma_{k_q})=
\sigma_{k_1}\cdots\wh{\sigma}_{k_i}\cdots\sigma_{k_q}=\sigma_{k_1}\cdots\sigma_{k_{i-1}}
\sigma_{k_{i+1}}\cdots\sigma_{k_q}.
\end{equation*}
Notice from Definition 2.4 that, the total subset $S_\uu=S_{\partial_i(\uu)}\cup (S_{\uu}\setminus S_{\partial_i(\uu)})$,
we define a $\kk[\xx]$-module homomorphism
$$\w{d}: \Lambda^{q,*}[\PP]\otimes_\kk \kk[\xx]\longrightarrow \Lambda^{q-1,*}[\PP]\otimes_\kk \kk[\xx]$$
by setting
\begin{equation}\tag{\nr}
\w{d}(\uu\otimes\xx^\ba)=\sum_{i=1}^q(-1)^i\partial_i(\uu)\otimes
    \xx_{(S_\uu\setminus S_{\partial_i(\uu)})}\xx^\ba,
\end{equation}
where $\displaystyle{\xx_{(S_\uu\setminus S_{\partial_i(\uu)})}=\frac{\xx_{S_\uu}}{\xx_{S_{\partial_i(\uu)}}}}$
is the monomial of $\kk[\xx]$ corresponding to the subset
$S_\uu\setminus S_{\partial_i(\uu)}$.

{\bf Theorem \nr \  (Taylor resolution)} {\it
$(\Lambda^{*,*}[\PP]\otimes_\kk\kk[\xx],\ \ \w{d})$ is a free resolution of $\kk[\xx]/\II_\PP$.
Thus
\[
H_{q,\ba}\left((\Lambda^{*,*}[\PP]\otimes_\kk\kk[\xx])\otimes_{\kk[\xx]}\kk,\ \  \w{d}\otimes 1\right)=
  Tor^{\kk[\xx]}_{q, \ba}(\kk[\xx]/\II_\PP, \kk).
\]
}
\begin{proof}  One may find a proof of this Theorem in Chapter 4 of \cite{Mil}, where it is given by the cellular resolution.
Here we give a proof by induction on the cardinal number $|\PP|=s$.

    Notice that, $\partial_i\partial_j=\partial_{j-1}\partial_i$ for $i<j$, an standard argument
shows that $\w{d}$ is a derivation, that is  $\w{d}\cdot\w{d}=0$.

    For $\Lambda^{0,*}[\PP]\otimes_\kk\kk[\xx]=\kk[\xx]$, we define
the augmentation
\[
\varepsilon: \Lambda^{0,*}[\PP]\otimes_\kk\kk[\xx]\longrightarrow \kk[\xx]/\II_\PP
\]
to be the quotient map. Furthermore from $\w{d}(\sigma_i)=-\xx_{\sigma_i}$
we see that
\[
\xymatrix {\Lambda^{1,*}[\PP]\otimes_\kk\kk[\xx] \ar[r]^-{\w{d}} & \Lambda^{0,*}[\PP]\otimes_\kk\kk[\xx]
  \ar[r]^-{\varepsilon} & \kk[\xx]/\II_\PP\ar[r] & 0}
\]
is exact. Thus we only need to prove that $H_{q,*}(\Lambda^{*,*}[\PP]\otimes_\kk\kk[\xx], \w{d})=0$ for all $q>0$,
and this will be done by induction on the cardinal number $|\PP|=s$.

  If $|\PP|=1$, that is $\PP=\{\sigma_1\}$, it is easy to see that
\[
\xymatrix @C=1.0cm @R=0.4cm
{
0 \ar[r] & \kk\{\sigma_1\}\otimes_\kk\kk[\xx] \ar[r]^-{\w{d}} \ar@{=}[d]
 & \kk[\xx] \ar[r]^-{\varepsilon} \ar@{=}[d] & \kk[\xx]/\II_\PP \\
0 \ar[r] & \Lambda^{1,*}[\{\sigma_1\}]\otimes_\kk\kk[\xx] \ar[r]^-{\w{d}} &
  \Lambda^{0,*}[\{\sigma_1\}]\otimes_\kk\kk[\xx] \ar[r]^-{\varepsilon} &
  \kk[\xx]/\II_\PP
}
\]
is exact and then $H_{1,*}(\Lambda^{*,*}[\{\sigma_1\}]\otimes_\kk\kk[\xx],\ \w{d})=0$.

    Inductively suppose that for any simplicial complement $\PP$ with cardinal number $|\PP|=s$, that is
$\PP=\{\sigma_1,\sigma_2,\cdots,\sigma_s\}$ contains $s$ subsets,
the homology $H_{q,*}(\Lambda^{*,*}[\PP]\otimes_\kk \kk[\xx], \w{d})=0$ for $q>0$.
Then for a simplicial complement $\mathbb{Q}$ with $|\mathbb{Q}|=s+1$,
that is
\[
\mathbb{Q}=\{\sigma_1,\sigma_2,\cdots,\sigma_s,\sigma\}=\PP\cup\{\sigma\}.
\]
$\Lambda^{*,*}[\PP]\otimes_\kk\kk[\xx]$ is a sub-complex of $\Lambda^{*,*}[\mathbb{Q}]\otimes_\kk\kk[\xx]$,
then the short exact sequence
\[
\xymatrix {0\ar[r] & \Lambda[\PP]\otimes_\kk\kk[\xx] \ar[r]^-{i} &
  \Lambda[\mathbb{Q}]\otimes_\kk\kk[\xx] \ar[r]
  & (\Lambda[\mathbb{Q}]/\Lambda[\PP])\otimes_\kk\kk[\xx]\ar[r] & 0}
\]
induces a long exact sequence in homologies
\[
\xymatrix @C=0.5cm
{\cdots\ar[r] & H_{q,*}\left(\Lambda[\PP]\otimes_\kk\kk[\xx]\right) \ar[r]^-{i_*} &
  H_{q,*}\left(\Lambda[\mathbb{Q}]\otimes_\kk\kk[\xx]\right) \ar[r]
  & H_{q,*}\left((\Lambda[\mathbb{Q}]/\Lambda[\PP])\otimes_\kk\kk[\xx]\right) \ar[r]^-{\partial} & \cdots}
\]
The quotient complex $(\Lambda^{q+1,*}[\mathbb{Q}]/\Lambda^{q+1,*}[\PP])\otimes_\kk\kk[\xx]$ is expressed as
$\Lambda^{q,*}[\PP]\sigma\otimes_\kk\kk[\xx]$ with generators $\uu\cdot\sigma\otimes\xx^{\ba}$ and  the  induced
differential
\[
\wh{d}: \Lambda^{q,*}[\PP]\sigma\otimes_\kk\kk[\xx]\rightarrow \Lambda^{q-1,*}[\PP]\sigma\otimes_\kk\kk[\xx]
\]
is given by
\[
\wh{d}(\uu\cdot\sigma\otimes\xx^\ba)=\sum_{i=1}^q(-1)^i(\partial_i(\uu))\cdot\sigma\otimes
  \xx_{S_i}\xx^\ba
\]
where $S_i=S_{\uu\cdot\sigma}\setminus S_{(\partial_i(\uu))\cdot\sigma}=
(\sigma_{k_1}\cup\cdots\cup\sigma_{k_q}\cup\sigma)\setminus(\sigma_{k_1}\cup\cdots
 \cup\wh{\sigma}_{k_i}\cup\cdots\cup\sigma_{k_q}\cup\sigma)$.

   By the induction hypothesis, we know that $H_{q,*}(\Lambda[\PP]\otimes_\kk\kk[\xx], \w{d})=0$ for $q>0$. Thus we
only need to prove that $H_{q,*}\left((\Lambda[\mathbb{Q}]/\Lambda[\PP])\otimes_\kk\kk[\xx]\right)=0$ for $q>1$ and
\[
\xymatrix @C=0.6cm
 {0\ar[r] & H_{1,*}\left((\Lambda[\mathbb{Q}]/\Lambda[\PP]\otimes_\kk\kk[\xx])\right) \ar[r]^-{\partial} &
  H_{0,*}(\Lambda[\PP]\otimes_\kk\kk[\xx])=\kk[\xx]/\II_\PP
 }
\]
is exact or equivalently $\partial$ is a monomorphism.  To do so,
let
\[
\PP'=\{\sigma'_1=\sigma_1\setminus\sigma, \sigma'_2=\sigma_2\setminus\sigma,\cdots,\sigma'_s=\sigma_s\setminus\sigma\}.
\]
Define a homomorphism of $\kk[\xx]$-modules
$f:\Sigma\Lambda[\PP']\otimes_\kk\kk[\xx]\rightarrow (\Lambda[\mathbb{Q}]/\Lambda[\PP])\otimes_\kk\kk[\xx]$
by
\[
f(\sigma'_{k_1}\sigma'_{k_2}\cdots\sigma'_{k_q})=\sigma_{k_1}\sigma_{k_2}\cdots\sigma_{k_q}\sigma.
\]
Notice that,
\begin{align*}
  & \left(\sigma_{k_1}\cup\cdots\cup\sigma_{k_q}\cup\sigma\right)\setminus
    \left(\sigma_{k_1}\cup\cdots\cup\wh{\sigma}_{k_i}\cup\cdots\cup\sigma_{k_q}\cup\sigma\right) \\
= & \left((\sigma_{k_1}\setminus\sigma)\cup\cdots\cup(\sigma_{k_q}\setminus\sigma)\right)
     \setminus\left((\sigma_{k_1}\setminus\sigma)\cup\cdots\cup\wh{(\sigma_{k_i}\setminus\sigma)}
     \cup\cdots\cup(\sigma_{k_q}\setminus\sigma)\right)\\
= & (\sigma'_{k_1}\cup\cdots\cup\sigma'_{k_q})\setminus(\sigma'_{k_1}
  \cup\cdots\cup\wh{\sigma}'_{k_i}\cup\cdots\cup\sigma'_{k_q})
\end{align*}
we see that for a monomial $\uu'=\sigma'_{k_1}\sigma'_{k_2}\cdots\sigma'_{k_n}\in \Sigma\Lambda^{*,*}[\PP']$,
$S_{\uu'}\setminus S_{\partial_i(\uu')}=S_{\uu\cdot\sigma}\setminus S_{(\partial_i(\uu))\cdot\sigma}$. Thus
\[
f:\Sigma\Lambda[\PP']\otimes_\kk\kk[\xx]\rightarrow \Lambda[\mathbb{Q}]/\Lambda[\PP]\otimes_\kk\kk[\xx]
\]
is an isomorphism of chain complexes. By induction on $s$ we see that
\[
H_{q+1,*}\left((\Lambda[\mathbb{Q}]/\Lambda[\PP])\otimes_\kk\kk[\xx],\ \wh{d}\right)
  \cong H_{q,*}(\Lambda[\PP']\otimes_\kk\kk[\xx],\ \w{d})=0
\]
for $q>0$ and
\[
H_{1,*}\left((\Lambda[\mathbb{Q}]/\Lambda[\PP])\otimes_\kk\kk[\xx],\ \wh{d}\right)
 \cong H_{0,*}(\Lambda[\PP']\otimes_\kk\kk[\xx],\ \w{d})
 \cong\kk[\xx]/\II_{\PP'}.
\]
$H_{0,*}(\Lambda[\PP']\otimes_\kk\kk[\xx],\ \w{d})\cong\kk[\xx]/\II_{\PP'}$ is
generated by $[\sigma\otimes \xx^\ba]$ with  $\xx^\ba\in\kk[\xx]$
and relations $[\sigma\otimes \xx^\ba]=0$ if $\xx^\ba\in \II_{\PP'}$.
The connecting homomorphism
\[
\xymatrix @C=0.6cm
 {\partial:
  H_{1,*}\left((\Lambda[\mathbb{Q}]/\Lambda[\PP])\otimes_\kk\kk[\xx],\ \wh{d}\right) \ar[r] &
  H_{0,*}\left(\Lambda[\PP]\otimes_\kk\kk[\xx], \w{d}\right)=\kk[\xx]/\II_\PP}
\]
is given by $\partial([\sigma\otimes\xx^\ba])=-[\xx_\sigma\xx^\ba]\in\kk[\xx]/\II_\PP$.

    Fix an $\NN^m$ vector $\mathbf{c}$, we see that
$H_{1,\mathbf{c}}\left((\Lambda[\mathbb{Q}]/\Lambda[\PP])\otimes_\kk\kk[\xx],\ \wh{d}\right)$
is generated by $[\sigma\otimes \xx^\ba]$ with an unique $\NN^m$ vector $\ba$ such that
$\xx_\sigma\xx^\ba=\xx^\mathbf{c}$.
If
\[
\partial([\sigma\otimes\xx^\ba])=-[\xx_\sigma\xx^\ba]=0\in
H_{0,*}(\Lambda[\PP]\otimes_\kk\kk[\xx],\ \w{d})=\kk[\xx]/\II_\PP,
\]
then in the $\mathbb{N}^m$ graded $\kk$-module $\kk[\xx]$, the monomial
$\xx_\sigma\xx^\ba\in\II_\PP$. This implies that there is a $\sigma_k\in \PP$ such that
$\xx_\sigma\xx^\ba=\xx_{\sigma_k}\xx^\mathbf{b}$. Thus $\xx_\sigma\xx^\ba$ could be expressed as
$\xx_{\sigma\cup\sigma_k}\xx^\mathbf{c}$.
Noticed that  $\sigma\cup\sigma_k=\sigma\cup\sigma'_k$ where $\sigma'_k=\sigma_k\setminus\sigma$,
one has $\xx_\sigma\xx^\ba=\xx_{\sigma}\xx_{\sigma'_k}\xx^\mathbf{c}$
and
\[
\wh{d}(\sigma'_k\sigma\otimes\xx^\mathbf{c})=-\sigma\otimes\xx_{\sigma'_k}\xx^\mathbf{c}=-\sigma\otimes\xx^\ba.
\]
This implies that $[\sigma\otimes\xx^\ba]=0\in H_{1,*}\left((\Lambda[\mathbb{Q}]/\Lambda[\PP])
\otimes\kk[\xx], \wh{d}\right)$
and the connecting homomorphism $\partial$ is a monomorphism.  The theorem follows.
\end{proof}

    To describe the differential
\[
\xymatrix @R=0.4cm{
  d:  & \Lambda^{q,*}[\PP]\otimes_\kk\kk[\xx]\otimes_{\kk[\xx]}\kk \ar[r]\ar@{=}[d]
      & \Lambda^{q-1,*}[\PP]\otimes_\kk\kk[\xx]\otimes_{\kk[\xx]}\kk  \ar@{=}[d] \\
      & \Lambda^{q,*}[\PP] \ar[r]
      & \Lambda^{q-1,*}[\PP], }
\]
let $\sigma$ be a subset of $[m]$ and let $\Lambda^{*,\sigma}[\PP]$ be the
submodules of $\Lambda[\PP]$ generated by the monomials $\uu$ with $S_\uu=\sigma$.
Then by Definition 2.4, $\Lambda[\PP]$ becomes an  $\NN\times 2^{[m]}$ graded module over $\kk$ and
$$\Lambda^{*,*}[\PP]=\bigoplus_{i,\sigma \in 2^{[m]}}\Lambda^{i,\sigma}[\PP].$$

    {\bf Theorem \nr} {\it Let $\PP=\{\sigma_1, \sigma_2, \cdots, \sigma_s\}$ be a simplicial complement
and $K_\PP$ be the corresponding simplecial complex by Definition 2.2.
Let $(\Lambda^{*,*}[\PP], d)$ be the chain complex induced from the bi-graded exterior algebra on $\PP$,
the differential
$d: \Lambda^{q,*}[\PP]\rightarrow \Lambda^{q-1,*}[\PP]$ is given by
\begin{align*}
d(\uu)= & \sum_i(-1)^i\partial_i(\uu)\cdot\delta_{\partial_i(\uu)}
\end{align*}
where $\delta_{\partial_i(\uu)}=1$ if the total subset $S_\uu=S_{\partial_i(\uu)}$
and $\delta_{\partial_i(\uu)}=0$ if $S_\uu\not=S_{\partial_i(\uu)}$. Then the homology of
$(\Lambda^{*,*}[\PP],\ d)$ is
\[
H_{*,\sigma}(\Lambda^{*,*}[\PP],\ d)=Tor^{\kk[\xx]}_{*,\sigma}(\kk[\xx]/\II_\PP, \kk)
  =Tor^{\kk[\xx]}_{*,\sigma}(\kk(K_\PP), \kk)
\]
which is called the homology of the simplicial complement $\PP$.}
\begin{proof}
    Consider the chain complex $(\Lambda^{q,*}[\PP]\otimes_\kk\kk[\xx]\otimes_{\kk[\xx]}\kk, \w{d}\otimes 1)$.
From (2.5), we see that for a generator $\uu\otimes 1$ of $\Lambda^{q,*}[\PP]\otimes_\kk\kk[\xx]$
\[
\w{d}(\uu\otimes 1)=\sum_{i=1}^q(-1)^i\partial_i(\uu)\otimes\xx_{(S_\uu\setminus S_{\partial_i(\uu)})}
\]
and $\xx_{(S_\uu\setminus S_{\partial_i(\uu)})}=1$ if and only if $S_\uu=S_{\partial_i(\uu)}$.
Notice that, $\kk[\xx]$ acts on $\kk$ by sending all $x_i$'s to $0$, the theorem follows from
\[
d(\uu\otimes_\kk 1\otimes_{\kk[\xx]} 1)=\w{d}(\uu\otimes_\kk 1)\otimes_{\kk[\xx]} 1=
 \sum_{i=1}^q (-1)^i\partial_i(\uu)\otimes_\kk \xx_{(S_\uu\setminus S_{\partial_i(\uu)})}\otimes_{\kk[\xx]}1
\]
and
\[
\partial_i(\uu)\otimes_\kk\xx_{(S_\uu\setminus S_{\partial_i(\uu)})}\otimes_{\kk[\xx]}1=0\in
 \Lambda^{q-1,*}[\PP]\otimes_\kk\kk[\xx]\otimes_{\kk[\xx]}\kk
\]
if $\xx_{S_\uu\setminus S_{\partial_i(\uu)}}\not=1$, that is $S_{\partial_i(\uu)}\not=S_\uu$.
\end{proof}

    Fix a $\sigma\in 2^{[m]}$ and let $\Lambda^{*,\sigma}[\PP]$ be the submodule
generated by the elements $\uu$ with total subset $S_\uu=\sigma$. From the definition
of $\delta_{\partial_i(\uu)}$ we see that $\partial_i(\uu)\cdot\delta_{\partial_i(\uu)}
\in \Lambda^{*,\sigma}[\PP]$, otherwise $\delta_{\partial_i(\uu)}=0$. Then $(\Lambda^{*,\sigma}[\PP],\ d)$
is a sub-complex of $\Lambda^{*,*}[\PP]$ and
\[
Tor_{q,\sigma}^{\kk[\xx]}(\kk[\xx]/\II_\PP,\ \kk)=H_q(\Lambda^{*,\sigma}[\PP], d)
\]

   A precise expression of the Hochster theorem (\cf \cite{BP}) is:

{\bf Theorem:} (Baskakov) {\it There are isomorphisms
$$Tor^{\kk[\xx]}_{q,\sigma}(\kk,\ \kk[\xx]/\II_{K_\PP})\cong \w{H}^{|\sigma|-q-1}(K_\PP\cap\sigma,\ \kk),$$
where $K_\PP\cap\sigma$ is the {\it full sub-complex} of $K_\PP$ consisting
of all simplices of $K_\PP$ which have all of their vertices in $\sigma$, that is
$K_\PP\cap\sigma=\{\tau\cap\sigma|\tau\in K_\PP\}$.}

    Thus we have the following theorem:

    {\bf Theorem \nr} (Combinatorial Hochster theorem) {\it Let $\PP$ and $K_\PP$ be as above, then
\begin{align*}
H_{q,\sigma}(\Lambda^{*,\sigma}[\PP],\ d)\cong & \w{H}^{|\sigma|-q-1}(K_\PP\cap\sigma,\ \kk)
& \mbox{and} & & H_{q,[m]}(\Lambda^{*,[m]}[\PP],\ d)\cong & \w{H}^{m-q-1}(K_\PP,\ \kk).
\end{align*}
Furthermore the cohomology module of the classical moment-angle complex is
\[
H^{2j-q}(\ZZ_{K_\PP}, \kk)=\bigoplus_{\sigma\subset [m], |\sigma|=j}
 H_{q,\sigma}(\Lambda^{*,*}[\PP],\ d).
\]}

  {\bf Remark:} Here we use the agreement that the cohomology of the empty simplicial complex $\{\phi\}$ is
$\w{H}^{-1}(\{\phi\}, \kk)=\kk$.
\begin{proof} From Theorem 2.7, we see that $H_{q,\sigma}(\Lambda^{*,\sigma}[\PP],\ d)=
Tor_{q,\sigma}^{\kk[\xx]}(\kk(K_\PP), \kk)$.  Then from
\[
Tor_{q,\sigma}^{\kk[\xx]}(\kk(K_\PP), \kk)=Tor_{q,\sigma}^{\kk[\xx]}(\kk,\ \kk(K_\PP))
  = \w{H}^{|\sigma|-q-1}(K_\PP\cap\sigma,\ \kk)
\]
we get the first one. The second one follows from the Hochster Theorem (\cf \cite{BP} for example).
\end{proof}

\section{The algebraic structure on $Tor_{i,\sigma}^{\kk[\xx]}(\kk(K_\PP), \kk)$}

    Let $A$ and $B$ be $\kk[\xx]$ algebras with structure maps
\begin{align*}
\varphi: & A\otimes_{\kk[\xx]}A\rightarrow A, & \varphi': B\otimes_{\kk[\xx]}B\rightarrow B.
\end{align*}
Choose a free resolution $P$ of $A$,
\[
\xymatrix @C=0.5cm{
P: \cdots \ar[r] & P_q \ar[r] & P_{q-1} \ar[r] & \cdots \ar[r] &
  P_1 \ar[r] & P_0 \ar[r] & A,
}
\]
$P\otimes P$ is a free resolution of $A\otimes_{\kk[\xx]}A$.
The structure map $\varphi:A\otimes_{\kk[\xx]}A$ gives a chain map
$\varphi: P\otimes_{\kk[\xx]} P\longrightarrow P$ which is unique up to chain homotopy.
\[
\xymatrix @C=0.4cm {
P\otimes_{\kk[\xx]} P: \ar[d]^-{\varphi} & \cdots \ar[r] & (P\otimes P)_q \ar[r] \ar[d]^-{\varphi_q} &
  (P\otimes P)_{q-1} \ar[r] \ar[d]^-{\varphi_{q-1}} & \cdots \ar[r] &
  (P\otimes P)_1 \ar[r] \ar[d]^-{\varphi_1} & (P\otimes P)_0 \ar[r] \ar[d]^-{\varphi_0} &
  A\otimes_{\kk[\xx]}  A \ar[d]^-{\varphi}\\
P:                                      & \cdots \ar[r] & P_q \ar[r] & P_{q-1} \ar[r] &
  \cdots \ar[r] & P_1 \ar[r] & P_0 \ar[r] & A
}
\]
Applying the functor $\otimes_{\kk[\xx]}(B\otimes_{\kk[\xx]}B)$, we get a natural chain map
\begin{equation}\tag{\nr}
\xymatrix @C=0.5cm{
(P\otimes_{\kk[\xx]} P)\otimes_{\kk[\xx]}(B\otimes_{\kk[\xx]}B) \ar@{=}[r] \ar[d]^-{\varphi\otimes 1} &
 (P\otimes_{\kk[\xx]}B)\otimes_{\kk[\xx]}(P\otimes_{\kk[\xx]}B) \ar[d]^-{\times} \\
P\otimes_{\kk[\xx]}(B\otimes_{\kk[\xx]}B) \ar[r]^-{1\otimes \varphi'} & P\otimes_{\kk[\xx]}B.
}
\end{equation}
On the other hand, there are natural maps
\begin{equation}\tag{\nr}
H_i\left(P\otimes_{\kk[\xx]}B\right)\otimes H_j\left(P\otimes_{\kk[\xx]}B\right)
 \longrightarrow H_{i+j}\left((P\otimes_{\kk[\xx]}B)\otimes_{\kk[\xx]}(P\otimes_{\kk[\xx]}B)\right);
\end{equation}
one maps the cycles $c\in (P\otimes_{\kk[\xx]}B)_i$ and $c'\in(P\otimes_{\kk[\xx]}B)_j$
to
\[
c\otimes c'\in \left((P\otimes_{\kk[\xx]}B)\otimes_{\kk[\xx]}(P\otimes_{\kk[\xx]}B)\right)_{i+j}.
\]
Applying the natural chain map $\times: (P\otimes_{\kk[\xx]}B)\otimes_{\kk[\xx]}(P\otimes_{\kk[\xx]}B)
\rightarrow P\otimes_{\kk[\xx]}B$, we get the algebraic structure
\begin{equation}\tag{\nr}
\xymatrix {
Tor_{i,*}^{\kk[\xx]}(A, B)\otimes Tor_{j,*}^{\kk[\xx]}(A,B) \ar[r]^-{\times} \ar@{=}[d] &
Tor_{i+j,*}^{\kk[\xx]}(A,B) \ar@{=}[d] \\
H_i(P\otimes_{\kk[\xx]}B)\otimes H_j(P\otimes_{\kk[\xx]}B) \ar[r] &
H_{i+j}(P\otimes_{\kk[\xx]}B)
}
\end{equation}

    Notice that both $\kk(K_\PP)$ and $\kk$ are $\kk[\xx]$ algebras,
$Tor_{i,*}^{\kk[\xx]}(\kk(K_\PP), \kk)$ is an algebra in a natural way and the
following isomorphism of algebras holds:
\[
H^*(\ZZ_{K_\PP}, \kk)\cong Tor_{i,*}^{\kk[\xx]}(\kk(K_\PP), \kk).
\]

    To describe the algebraic structure of $Tor_{i,*}^{\kk[\xx]}(\kk(K_\PP), \kk)$,
consider the free $\kk[\xx]$-resolution of $\kk(K_\PP)$ given by Theorem 2.6,
which is denoted by $R^{*,*}$
$$R^{*,*}=\Lambda^{*,*}[\PP]\otimes_\kk \kk[\xx].$$
The tensor product
\[
R^{*,*}\otimes_{\kk[\xx]}R^{*,*}\cong (\Lambda^{*,*}[\PP]\otimes_\kk\Lambda^{*,*}[\PP])\otimes_\kk
  (\kk[\xx]\otimes_{\kk[\xx]}\kk[\xx])\cong
  (\Lambda^{*,*}[\PP]\otimes_\kk\Lambda^{*,*}[\PP])\otimes_\kk\kk[\xx]
\]
gives a free resolution of $\kk(K_\PP)\otimes_{\kk[\xx]}\kk(K_\PP)\cong \kk(K_\PP)$. The algebraic structure
of $\kk(K_\PP)$ is given by
\[
\xymatrix {
\kk(K_\PP)\otimes_{\kk[\xx]}\kk(K_\PP)\cong \kk(K_\PP) \ar[r]^-{1} & \kk(K_\PP).
}
\]
Thus we need to construct a chain map
$\varphi: (\Lambda^{*,*}[\PP]\otimes_\kk\Lambda^{*,*}[\PP])\otimes_\kk \kk[\xx]\rightarrow
  \Lambda^{*,*}[\PP]\otimes_\kk \kk[\xx]$
that makes the following diagram commute
\[
\xymatrix {
(\Lambda^{*,*}[\PP]\otimes_\kk\Lambda^{*,*}[\PP])\otimes_\kk \kk[\xx] \ar[r]^-{\varphi} \ar[d]^-{\varepsilon} &
  \Lambda^{*,*}[\PP] \ar[d]^-{\varepsilon}\\
  \kk(K_\PP)\otimes_{\kk[\xx]} \kk(K_\PP) =\kk(K_\PP) \ar[r]^-{1} & \kk(K_\PP)
}
\]

{\bf Construction \nr} {\it The exterior algebra $\Lambda^{*,*}[\PP]$ has natural product structure.
Given two generators $\uu$ and $\vv$ of $\Lambda^{*,*}[\PP]$ the total subset of
$\uu\cdot\vv$ is $S_{\uu\cdot\vv}=S_{\uu}\cup S_{\vv}$ if $\uu\cdot\vv\not=0$ in
$\Lambda^{*,*}[\PP]$. Thus the monomial $\xx_{S_{\uu\cdot\vv}}$ is a factor of
$\xx_{S_\uu}\cdot\xx_{S_\vv}$ in $\kk[\xx]$.  We define the product
\[
\varphi: (\Lambda^{*,*}[\PP]\otimes_\kk\Lambda^{*,*}[\PP])\otimes_\kk\kk[\xx]
\longrightarrow \Lambda^{*,*}[\PP]\otimes_\kk\kk[\xx]
\]
to be the $\kk[\xx]$-module map induced by
\begin{align}\tag{\nr}
\uu\times\vv = & \uu\cdot\vv \otimes\frac{\xx_{S_\uu}\cdot\xx_{S_\vv}}{\xx_{S_{\uu\cdot\vv}}}.
\end{align}
It is apparent that so defined product keeps the second degree ($\NN^m$ degree).
}

{\bf Theorem \nr} {\it The product defined in Construction 3.4 is a chain map. Thus the algebraic structure
in $H_{*,*}(\Lambda^{*,*}[\PP],d)\cong Tor_{*,*}^{\kk[\xx]}(\kk(K_\PP), \kk)$ is given
by
\[
[c]\times [c']=[c\cdot c']\cdot \delta_{\sigma,\sigma'},
\]
where $[c]\in H_{q,\sigma}(\Lambda^{*,\sigma}[\PP], d)$, $[c']\in H_{j,\sigma'}(\Lambda^{*,\sigma'}[\PP], d)$,
the cycle $c\cdot c'$ is product of $c$ and $c'$ in $\Lambda^{*,*}[\PP]$ and $ \delta_{\sigma,\sigma'}=1$
if $\sigma\cap\sigma'=\phi$, $ \delta_{\sigma,\sigma'}=0$ if $\sigma\cap\sigma'\not=\phi$.
}
\begin{proof}
Notice that $\Lambda^{0,*}[\PP]=\kk$, it is apparent that
\[
\xymatrix{
(\Lambda^{0,*}[\PP]\otimes_\kk\Lambda^{0,*}[\PP])\otimes_\kk \kk[\xx]
\ar^-{\varepsilon}[r]\ar^-{\varphi_0}[d] &
\kk(K_\PP)\otimes_{\kk[\xx]} \kk(K_\PP)=\kk(K_\PP) \ar^-{1}[d]\\
\Lambda^{0,*}[\PP]\otimes_\kk \kk[\xx] \ar^-{\varepsilon}[r] & \kk(K_\PP)}
\]
commute.

    To prove that $\varphi$ defined by Construction 3.4 is a chain map, we only need to check that
\[
d(\uu\times\vv)=d(\uu)\times\vv+(-1)^q\uu\times d(\vv)
\]
or equivalently to check that $\Lambda^{*,*}[\PP]\otimes_\kk \kk[\xx]$ is a differential graded algebra,
because the differential
in $(\Lambda^{*,*}[\PP]\otimes_\kk\Lambda^{*,*}[\PP])\otimes_\kk \kk[\xx]$ is given by
$d(\uu\otimes\vv)=d(\uu)\otimes\vv+(-1)^q\uu\otimes d(\vv)$ if $\uu\in\Lambda^{q,*}[\PP]$.

    Let $V$ be the free $\kk$-module generated by $\PP=\{\sigma_1, \sigma_2, \cdots, \sigma_s\}$
and let $T(V)$ be the tensor algebra of $V$, that is
$\displaystyle{T(V)=\bigoplus_{q\geqslant 0} V^{\otimes q}}$. We still use
\begin{align*}
\uu= & \sigma_{k_1}\sigma_{k_2}\cdots\sigma_{k_q} & &\mbox{and} &
S_\uu= & \sigma_{k_1}\cup\sigma_{k_2}\cup\cdots\cup\sigma_{k_q}
\end{align*}
to denote the generator of $T(V)$ and the total subset of $\uu$ respectively. The exterior algebra
$$\Lambda^{*,*}[\PP]=T(V)/\langle \uu\vv-(-1)^{pq}\vv\uu, \uu\uu\rangle$$
where $\langle \uu\vv-(-1)^{pq}\vv\uu, \uu\uu\rangle$ is the ideal generated by
$\uu\vv-(-1)^{pq}\vv\uu$ and $\uu\uu$ for any generators $\uu\in T^{q, *}(V)$ and
$\vv\in T^{p,*}(V)$.

    Similarly, the monomial $\xx_{S_\uu}\in \kk[\xx]$ corresponding to the total subset
$S_\uu$ is expressed as $\xx_{S_\uu}=\xx^{\ba_\uu}$ for a unique $\NN^m$ vector $\ba_\uu$.
We also define
\begin{align*}
bideg(\uu) = & (q, \ba_\uu) & & \mbox{and} \\
\partial_i(\uu) = & \sigma_{k_1}\cdots \wh{\sigma}_{k_i}\cdots \sigma_{k_q}.
\end{align*}
Consider the $\NN\times\NN^m$ graded free $\kk[\xx]$-module $T(V)\otimes_\kk \kk[\xx]$,
we define a differential
$$\w{d}: T^{q,*}(V)\otimes_\kk \kk[\xx]\rightarrow
T^{q-1,*}(V)\otimes_\kk \kk[\xx]$$
by setting
$$\w{d}(\uu\otimes\xx^{\ba})=\sum_i(-1)^i\partial_i(\uu)\otimes \frac{\xx_{S_\uu}}{\xx_{S_{\partial_i(\uu)}}}
  \cdot\xx^\ba
$$
as (2.5).

    Similar to Construction 3.4, $T(V)\otimes_\kk \kk[\xx]$ has algebraic structure given by
\[
\uu\times\vv=\uu\cdot\vv \otimes \frac{\xx_{S_\uu}\cdot\xx_{S_\vv}}{\xx_{S_{\uu\cdot\vv}}}.
\]
From
\begin{align*}
  & \w{d}(\uu\times\vv) = \w{d}(\uu\vv)\otimes \frac{\xx_{S_\uu}\cdot\xx_{S_\vv}}{\xx_{S_{\uu\vv}}}\\
= & \sum_i(-1)^i\partial_i(\uu)\vv\otimes\frac{\xx_{S_{\uu\vv}}}{\xx_{S_{\partial_i(\uu)\vv}}}\cdot
    \frac{\xx_{S_\uu}\cdot\xx_{S_\vv}}{\xx_{S_{\uu\vv}}}
+(-1)^q\sum_i(-1)^i\uu\partial_i(\vv)\otimes  \frac{\xx_{S_{\uu\vv}}}{\xx_{S_{\uu\partial_i(\vv)}}}\cdot
    \frac{\xx_{S_\uu}\cdot\xx_{S_\vv}}{\xx_{S_{\uu\vv}}}
\end{align*}
and
\begin{align*}
  & \w{d}(\uu)\times\vv + (-1)^q\uu\times\w{d}(\vv)\\
= & \sum_i(-1)^i\partial_i(\uu)\times\vv\otimes \frac{\xx_{S_\uu}}{\xx_{S_{\partial_i(\uu)}}}
+(-1)^q\sum_i(-1)^i\uu\times\partial_i(\vv)\otimes \frac{\xx_{S_\vv}}{\xx_{S_{\partial_i(\vv)}}}\\
= & \sum_i(-1)^i\partial_i(\uu)\vv\otimes
    \frac{\xx_{S_{\partial_i(\uu)}}\cdot\xx_{S_\vv}}{\xx_{S_{\partial_i(\uu)\vv}}}\cdot
    \frac{\xx_{S_\uu}}{\xx_{S_{\partial_i(\uu)}}}
+(-1)^q\sum_i(-1)^i\uu\partial_i(\vv)\otimes
    \frac{\xx_{S_\uu}\cdot\xx_{S_{\partial_i(\vv)}}}{\xx_{S_{\uu\partial_i(\vv)}}}\cdot
    \frac{\xx_{S_\vv}}{\xx_{S_{\partial_i(\vv)}}}
\end{align*}
we see that $\w{d}(\uu\times\vv)=\w{d}(\uu)\times\vv+(-1)^q\uu\times\w{d}(\vv)$ and
 $T(V)\otimes_\kk \kk[\xx]$ is a differential graded algebra.  It is easy to check
that the ideal $\langle \uu\vv-(-1)^{pq}\vv\uu, \uu\uu\rangle$ is invariant under the differential
$\w{d}$. Thus $\{T(V)\otimes_\kk \kk[\xx], \w{d}\}$ induces the differential graded algebraic structure
on $\Lambda^{*,*}[\PP]\otimes_\kk \kk[\xx]$ as desired.

    To prove the second part of the Theorem, applying the functor $\otimes_{\kk[\xx]}\kk$,
 we see that in $(\Lambda^{*,*}[\PP]\otimes_\kk \kk[\xx])\otimes_{\kk[\xx]}\kk$
\[
(\uu\times\vv\otimes_\kk 1)\otimes_{\kk[\xx]} 1=(\uu\cdot\vv\otimes_\kk 1)\otimes_{\kk[\xx]}
 \frac{\xx_{S_\uu}\cdot\xx_{S_\vv}}{\xx_{S_{\uu\vv}}}=0
\]
if the monomial $\displaystyle{\frac{\xx_{S_\uu}\cdot\xx_{S_\vv}}{\xx_{S_{\uu\vv}}}\not=1}$
or equivalently $S_\uu\cap S_\vv\not=\phi$.
\end{proof}

   {\bf Example:} Compute the $Tor^{\kk[\xx]}_{*,*}(\kk(K), \kk)$ and  the cohomology of $\ZZ_K$
corresponding to $K$ (See Figure 1)

  There is the simplicial complement $\PP=\left\{\sigma_1=\{1,5\},\ \sigma_2=\{2,4\},
\ \sigma_3=\{1,2,3\},\ \sigma_4=\{3,4,5\}\right\}$ corresponding to the simplicial complex  $K=K_\PP$.

\begin{center}
\setlength{\unitlength}{0.3mm}
\begin{picture}(210,120)(0,-10)
\put(0,0){\line(1,0){200}}
\put(0,0){\line(1,1){100}}
\put(200,0){\line(-1,1){100}}
\put(66.6,33.3){\line(1,0){66.6}}
\put(66.6,33.3){\line(1,2){33.3}}
\put(100,100){\line(1,-2){33.3}}
\put(0,0){\line(2,1){66.6}}
\put(200,0){\line(-2,1){66.6}}
\put(-7,-5){\mbox{1}}
\put(205,-5){\mbox{2}}
\put(66.6,20){\mbox{4}}
\put(130,20){\mbox{5}}
\put(98,103){\mbox{3}}
\put(80,-15){\mbox{Figure 1}}
\put(10,10){\line(0,-1){5}}
\put(20,20){\line(0,-1){10}}
\put(30,30){\line(0,-1){15}}
\put(40,40){\line(0,-1){20}}
\put(50,50){\line(0,-1){25}}
\put(60,60){\line(0,-1){30}}
\put(70,70){\line(0,-1){30}}
\put(80,80){\line(0,-1){20}}
\put(90,90){\line(0,-1){9}}
\put(110,90){\line(0,-1){9}}
\put(120,80){\line(0,-1){20}}
\put(130,70){\line(0,-1){30}}
\put(140,60){\line(0,-1){30}}
\put(150,50){\line(0,-1){25}}
\put(160,40){\line(0,-1){20}}
\put(170,30){\line(0,-1){15}}
\put(180,20){\line(0,-1){10}}
\put(190,10){\line(0,-1){5}}
\end{picture}
\end{center}\vspace{5mm}

From Theorem 2.7, we compute the differentials in $\Lambda^{*,*}[\PP]$;
\begin{align*}
d(\sigma_1\sigma_2\sigma_3\sigma_4)= & -\sigma_2\sigma_3\sigma_4 +\sigma_1\sigma_3\sigma_4
                                       -\sigma_1\sigma_2\sigma_4 +\sigma_1\sigma_2\sigma_3;\\
d(\sigma_1\sigma_3\sigma_4) = & -\sigma_3\sigma_4.
\end{align*}
And these are all the non-trivial differentials. Thus the $Tor^{\kk[\xx]}_{*,*}(\kk(K), \kk)$
is a free $\kk$-module generated by
\[
\left\{\begin{array}{ll}
         \sigma_1\sigma_2\sigma_4\\
         \sigma_1\sigma_2\sigma_3
         \end{array}
        \begin{array}{ll}
         \sigma_1\sigma_2\\
         \sigma_1\sigma_3\\
         \sigma_1\sigma_4\\
         \sigma_2\sigma_3\\
         \sigma_2\sigma_4
        \end{array}
        \begin{array}{ll}
         \sigma_1\\
         \sigma_2\\
         \sigma_3\\
         \sigma_4
        \end{array}
        \begin{array}{ll}
         1
        \end{array}
\right\}.
\]

    Consider the algebraic structure of $Tor_{*,*}^{\kk[\xx]}(\kk(K), \kk)$, we see from Theorem 3.6
that
\[
\sigma_1\times\sigma_2=\sigma_1\sigma_2
\]
is the only non-trivial product.  By the Hochester Theorem, we see that $H^*(\ZZ_K, \kk)$ is a free $\kk$-module with Poincar\`{e}
series $1+2x^3+2x^5+5x^6+2x^7$.

\section{The simplicial complements and the simplicial complexes}

   In \cite{CL} X. Cao and Z.~L\"{u} introduced a $Z/2$-algebra
$$2^{[m]^*}=\{f|f:2^{[m]}\rightarrow Z/2=\{0,1\}\}$$
consists of all $Z/2$ valued functions on the power set $2^{[m]}$ with addition $(f+g)(\sigma)=f(\sigma)+g(\sigma)$
and multiplication $(f\cdot g)(\sigma)=f(\sigma)\cdot g(\sigma)$.

    Let $f: 2^{[m]}\rightarrow Z/2$ be a function in $2^{[m]*}$.  We define
$$supp(f)=\{\sigma\in 2^{[m]}\mid f(\sigma)=1\}=f^{-1}(1)\subset 2^{[m]}$$
which is called the {\it support} of $f$.

    On the other hand, given a simplicial complex $K$,
there is the {\it characteristic function} $f_K\in 2^{[m]^*}$ of $K$ defined by
\begin{align}\tag{\nr}
f_K(\sigma)= & \left\{\begin{array}{ll}
                         1 & \mbox{if $\sigma\in K$}\\
                         0 & \mbox{otherwise}
                  \end{array}\right. ,
\end{align}
with $supp(f_K)=K$.

    Given a simplex $\sigma\in 2^{[m]}$, there are functions $\delta_\sigma$ and $\mu_\sigma\in 2^{[m]*}$
defined as
\begin{align}\tag{\nr}
\delta_\sigma(\tau)= & \left\{\begin{array}{ll}
                       1 & \mbox{if $\tau=\sigma$}\\
                       0 & \mbox{otherwise.}
                             \end{array}\right.
& \mu_{\sigma}(\tau)= & \left\{\begin{array}{ll}
                           1 & \mbox{if $\sigma\subset\tau$}\\
                           0 & \mbox{otherwise}
                          \end{array}\right.
\end{align}
From (4.2), one can easily show that any $f\in 2^{[m]*}$ can be expressed as
\[
f=\sum_{\sigma\in 2^{[m]}}f(\sigma)\delta_\sigma=\sum_{\sigma\in supp(f)}
  \delta_\sigma.
\]
Thus $\{\delta_\sigma\mid \sigma\in 2^{[m]}\}$ forms a basis of $2^{[m]^*}$.

{\bf Theorem \nr} {\it Let $\PP=\{\sigma_1,\sigma_2,\cdots,\sigma_s\}$ be a simplicial complement,
$K_\PP$ be the simplicial complex corresponding to $\PP$ by Definition 2.2 and $f_{K_\PP}$ be
the  characteristic function corresponding to
$K_\PP$ by (4.1). Then
\[
f_{K_\PP}=\prod_{i=1}^{s}(1+\mu_{\sigma_i})=\sum_{\uu\in\Lambda[\PP]}\mu_{S_\uu}
=\sum_{\sigma\in 2^{[m]}}\left(\sum_i dim_\kk\Lambda^{i,\sigma}[\PP]\right)\mu_\sigma,
\]
where the first sum runs over all the $\kk$-module generators $\uu$ of $\Lambda^{*,*}[\PP]$.}
\begin{proof}
    Notice that, $1(\tau)=1$ for any $\tau\in 2^{[m]}$ and $(1+\mu_{\sigma})(\tau)=0$
if and only if $\sigma\subset\tau$, we see that for any $\tau\in 2^{[m]}$
\[
f_{K_\PP}(\tau)=(1+\mu_{\sigma_1})(\tau)\cdot(1+\mu_{\sigma_2})(\tau)\cdot\cdots\cdot
       (1+\mu_{\sigma_s})(\tau)=0
\]
if and only if there exists an $i$ such that $(1+\mu_{\sigma_i})(\tau)=0$.
This is equivalent to that there is a $\sigma_i\in\PP$ such that $\sigma_i\subset\tau$.
By Definition 2.2, $\tau$ is a non-face of $K_\PP$. Thus
$\displaystyle{f_{K_\PP}=\prod_{i=1}^{s}(1+\mu_{\sigma_i})}$.

    For the second part of the theorem, notice that
for any subset $\omega\in 2^{[m]}$,
$(\mu_{\sigma}\cdot\mu_{\tau})(\omega)=\mu_{\sigma}(\omega)\cdot\mu_{\tau}(\omega)=1$
if and only if $\sigma\subset\omega$ and $\tau\subset\omega$ $\Longleftrightarrow$
$\sigma\cup\tau\subset\omega$. One has
\[
\mu_{\sigma}\cdot\mu_{\tau}=\mu_{\sigma\cup\tau}.
\]
Thus
\begin{align*}
\prod_{i=1}^{s}(1+\mu_{\sigma_i})= &
    1+\sum_{\{i_1,i_2\cdots,i_q\}\subset [m]}\mu_{\sigma_{i_1}}\mu_{\sigma_{i_2}}\cdots\mu_{\sigma_{i_q}}
    = \sum_{\uu\in\Lambda[\PP]}\mu_{S_\uu}
\end{align*}
by Definition 2.4.

    Fix a $\sigma\in 2^{[m]}$, $\Lambda^{*,\sigma}[\PP]$ is the $\kk-$submodule generated
by $\uu$ with total subset $S_\uu=\sigma$,
$$\sum_{S_\uu=\sigma}\mu_{S_\uu}=
 \left(\sum_i \mbox{dim}_\kk\Lambda^{i,\sigma}[\PP]\right)\mu_\sigma.$$
The third part of the Theorem follows.
\end{proof}

  {\bf Definition \nr} {\it Fix a subset $\omega\in [m]$, there is the map $\varepsilon_\omega: 2^{[m]}\longrightarrow 2^{[m]}$ defined
by $\varepsilon_\omega(\tau)=\tau\cup\omega$. The $\omega$-compression
$E_\omega: 2^{[m]^*}\longrightarrow 2^{[m]^*}$ is defined by
\[
E_\omega(f)=f\circ\varepsilon_\omega.
\]
For the simplicial complement $\PP=\{\sigma_1, \sigma_2, \cdots, \sigma_s\}$,  we
define the  $\omega$-compression of $\PP$ by
\[
E_\omega\PP=\{\sigma_1\setminus \omega, \sigma_2\setminus \omega, \cdots, \sigma_s\setminus \omega\}.
\]}

    From $(f+g)\circ\varepsilon_\omega=f\circ\varepsilon_\omega+g\circ\varepsilon_\omega$ and
$(f\cdot g)\circ\varepsilon_\omega=(f\circ\varepsilon_\omega)\cdot(g\circ\varepsilon_\omega)$,
we see that $E_\omega: 2^{[m]^*}\longrightarrow 2^{[m]^*}$ is a homomorphism of $Z/2$-algebra.

    {\bf Lemma \nr} {\it Let $\PP$ be a simplicial complement and $\displaystyle{f_{K_\PP}=\prod_{i=1}^s(1+\mu_{\sigma_i})}$
be the characteristic function of the simplicial complex $K_\PP$ by Theorem 4.3. Then
\[
E_\omega(f_{K_\PP})=\prod_{i=1}^s(1+\mu_{\sigma_i\setminus\omega})=f_{K_{E_\omega\PP}},
\]
i.e. $E_\omega(f_{K_\PP})$ is the characteristic function of the simplicial complex $K_{E_\omega\PP}$
corresponding to the simplicial complement $E_\omega\PP$.
}
\begin{proof}
    The $\omega$-compression $E_\omega: 2^{[m]^*} \longrightarrow 2^{[m]^*}$ is a homomorphism
of $Z/2$-algebra. Thus
\[
E_\omega(f_{K_\PP})=E_\omega\left(\prod_{i=1}^s(1+\mu_{\sigma_i})\right)=
    \prod_{i=1}^s\left(E_\omega(1)+E_\omega(\mu_{\sigma_i})\right).
\]

    For the function $\mu_\sigma$, take a $\tau\in 2^{[m]}$, we have
\[
E_\omega(\mu_\sigma)(\tau)=1\Longleftrightarrow\mu_{\sigma}(\varepsilon_\omega(\tau))=1
 \Longleftrightarrow\sigma\subset\tau\cup \omega\Longleftrightarrow(\sigma\setminus \omega)\subset\tau.
\]
This implies that
\begin{align}\tag{\nr}
E_\omega(\mu_{\sigma})= & \mu_{\sigma\setminus \omega}.
\end{align}
Then from
$E_\omega(1)=1$, we see that
$\displaystyle{E_\omega(f_{K_\PP})=\prod_{i=1}^s\left(1+E_\omega(\mu_{\sigma_i})\right)=
 \prod_{i=1}^s\left(1+\mu_{\sigma_i\setminus \omega}\right)=f_{K_{E_\omega\PP}}}.$
\end{proof}

    For an arbitrary simplex $\omega\in K$, define its {\it link} and {\it star} as the
sub-complexes
\begin{align*}
  & \mbox{star}_K\omega = \{\tau\in K|\omega\cup\tau\in K\};
  & \mbox{link}_K\omega = \{\tau\in K|\omega\cup\tau\in K, \omega\cap\tau=\phi\}.
\end{align*}

  {\bf Theorem \nr} {\it Let $\PP$ be a simplicial complement and $K_\PP$ be the simplicial complex corresponding
to $\PP$. Then:}
\begin{enumerate}
\item {\it The $star$,} $\mbox{star}_{K_\PP}\omega=K_{E_\omega\PP}$ {\it is the simplicial
complex corresponding to the simplicial complement $E_\omega\PP$. Thus}
\[
H_{q,\sigma}(\Lambda^{*,\sigma}[E_\omega\PP], d)=
  Tor^{\kk[\xx]}_{q,\sigma}(\kk(\mbox{star}_{K_\PP}\omega), \kk),
\]
\item {\it The $link$,} $\mbox{link}_{K_\PP}\omega=K_{E_\omega\PP}\cap ([m]\setminus\omega)$ {and}
\[
H_{q,[m]\setminus \omega}(\Lambda^{*,[m]\setminus \omega}[E_\omega\PP], d)=
 \w{H}^{m-|\omega|-q-1}(\mbox{link}_{K_\PP}\omega, \kk).
\]
\end{enumerate}

\begin{proof}
    Let $f_{K_\PP}$ and $f_{K_{E_\omega\PP}}$ be the characteristic functions of $K_\PP$ and
$K_{E_\omega\PP}$ respectively. Then a simplex $\tau\in \mbox{star}_{K_\PP}\omega$ if and only if
\begin{align*}
\tau\cup\omega\in K_\PP & \Longleftrightarrow f_{K_\PP}(\tau\cup\omega)=1
\Longleftrightarrow f_{K_\PP}\circ\varepsilon_\omega(\tau)=1
\Longleftrightarrow E_\omega(f_{K_\PP})(\tau)=f_{K_{E_\omega\PP}}(\tau)=1
\end{align*}
Thus the characteristic function of $\mbox{star}_{K_\PP}\omega$ is $f_{K_{E_\omega\PP}}$ and
\begin{align*}
\mbox{star}_{K_\PP}\omega= & \{\tau\in K_\PP|\omega\cup\tau\in K_\PP\}= K_{E_\omega\PP},  \\
\mbox{link}_{K_\PP}\omega= & \{\tau\in K_\PP|\omega\cup\tau\in K_\PP, \tau\cap\omega=\phi\} =
  K_{E_\omega\PP}\cap([m]\setminus\omega).
\end{align*}
From Theorem 2.7 we see that
\begin{align*}
Tor^{\kk[\xx]}_{q,\sigma}(\kk(\mbox{star}_{K_\PP}\omega), \kk)= & Tor^{\kk[\xx]}_{q,\sigma}(\kk(K_{E_\omega\PP}), \kk)
    =H_{q,\sigma}(\Lambda^{*,\sigma}[E_\omega\PP], d)
\end{align*}
and from Theorem 2.8 we get
\begin{align*}
\w{H}^{m-|\omega|-q-1}(\mbox{link}_{K_\PP}\omega,\ \kk) = &
\w{H}^{m-|\omega|-q-1}(K_{E_\omega\PP}\cap([m]\setminus\omega),\ \kk) =
H_{q,[m]\setminus \omega}(\Lambda^{*,[m]\setminus \omega}[E_\omega\PP], d).
\end{align*}
The theorem follows.
\end{proof}

\section{The cohomology of the generalized moment-angle complex}

    In this section we consider the cohomology module of the generalized moment-angle complex.
Recall from \cite{BB} \cite{BBCG08},

  {\bf Definition \nr} {\it Let $(\underline{X}, \underline{A})=\{(X_i, A_i, x_i)|i\in[m]\}$ denote
a set of triples of $CW$-complexes with base points $x_i\in A_i$ and $K_\PP$ be an abstract simplicial
complex. The generalized moment-angle complex determined by $(\underline{X}, \underline{A})$ and
$K_\PP$ denoted by $\ZZ_{K_\PP}(\underline{X}, \underline{A})$ is defined using the functor
\[
D: K_\PP\longrightarrow CW_*
\]
as follows: For every $\omega\in K_\PP$, let
\begin{align*}
D(\omega)= & Y_1\times Y_2 \times \cdots \times Y_m &\mbox{where} & &
Y_i=& \left\{\begin{array}{lll}
          X_i & \mbox{if $i\in\omega$}\\
          A_i & \mbox{if $i\not\in\omega$.}
            \end{array}\right.
\end{align*}
The generalized moment-angle complex is
\[
\ZZ_{K_\PP}(\underline{X}, \underline{A})=\bigcup_{\omega\in K_\PP} D(\omega).
\]}

    Let $Y_1\wedge Y_2 \wedge \cdots \wedge Y_m $ be the smash product given by
the quotient space
\[
Y_1\times Y_2\times\cdots\times Y_m/S(Y_1\times Y_2\times\cdots\times Y_m)
\]
where $S(Y_1\times Y_2\times\cdots\times Y_m)$ is the subspace of the product with at least
one coordinate given by the base-point $x_i\in Y_i$. The {\it generalized smash moment-angle
complex} is defined to be the image of $\ZZ_{K_\PP}(\underline{X}, \underline{A})$
in $X_1\wedge X_2\wedge\cdots\wedge X_m$, that is
\[
\wh{\ZZ}_{K_\PP}(\underline{X}, \underline{A})=\bigcup_{\omega\in K_\PP} \wh{D}(\omega)
\]
where
\begin{align*}
\wh{D}(\omega)= & Y_1\wedge Y_2 \wedge \cdots \wedge Y_m  &\mbox{subject to} & &
Y_i=& \left\{\begin{array}{lll}
          X_i & \mbox{if $i\in\omega$}\\
          A_i & \mbox{if $i\not\in\omega$.}
            \end{array}\right.
\end{align*}

    Given a non-empty subset $I=\{i_1, i_2, \cdots, i_k\}\subset [m]$ and a family of pairs
$(\underline{X}, \underline{A})$, define
\[
(\underline{X_I}, \underline{A_I})=\{(X_{i_j}, A_{i_j})|i_j\in I\}
\]
which is the subfamily of $(\underline{X}, \underline{A})$ determined by $I$.  It is known from
\cite{BB,BBCG08} that:
\begin{align}\tag{\nr}
H: \Sigma(\ZZ_{K_\PP}(\underline{X}, \underline{A}))\longrightarrow &
  \Sigma\left(\bigvee_{I\subset [m]}\wh{\ZZ}_{K_\PP\cap I}(\underline{X_I}, \underline{A_I})\right)
\end{align}
is a natural pointed homotopy equivalence.  To describe it more precisely, let
$\omega\subset I$ and
\begin{align*}
\wh{D}_I(\omega)= & Y_{i_1}\wedge Y_{i_2} \wedge \cdots\wedge Y_{i_k}
 & \mbox{with} & & Y_{i_j}=\left\{\begin{array}{ll}
                           X_{i_j} & \mbox{if $i_j\in \omega$}\\
                           A_{i_j} & \mbox{if $i_j\not\in\omega$.}
                   \end{array}\right.
\end{align*}
Then
\[
\wh{\ZZ}_{K_\PP\cap I}(\underline{X_I}, \underline{A_I})=\bigcup_{\omega\in K_\PP\cap I}\wh{D}_I(\omega).
\]

    Associated to a simplicial complex $K_\PP$, there is a partial ordered set ({\it poset}) $\overline{K_\PP}$
with point $\sigma$ in $\overline{K_\PP}$ corresponding to a simplex $\sigma\in K_\PP$ and order given by
reverse inclusion of simplices. Thus $\sigma_1\leqslant \sigma_2$ in $\overline{K_\PP}$ if and only if
$\sigma_2\subseteq\sigma_1$ in $K_\PP$.  Given an $\omega\in\overline{K_\PP}$ there are further posets given
by
\begin{align*}
\overline{K_\PP}_{<\omega} = & \{\tau\in\overline{K_\PP}|\tau<\omega\}=\{\tau\in\overline{K_\PP}|\omega\varsubsetneq\tau\}
 & \mbox{and} \\
\overline{K_\PP}_{\leqslant\omega} = & \{\tau\in\overline{K_\PP}|\tau\leqslant\omega\}
 =\{\tau\in\overline{K_\PP}|\omega\subset\tau\}.
\end{align*}

    Given a poset $P$, the {\it order complex} $\Delta(P)$ is the simplicial complex with vertices given by set of
points of $P$ and $k$-simplices given by ordered $(k+1)$-tuples $(p_0, p_1, \cdots, p_k)$ with
$p_0<p_1<\cdots<p_k$.  It follows that $\Delta(\overline{K_\PP})_{<\phi}=K'_{\PP}$ is the barycentric subdivision
of $K_\PP$.

    Given a simplicial complex $K_\PP$, we use $|K_\PP|$ to denote its geometric realization.
The symbol $|K_\PP|\ast Y$ denotes the join of $K_\PP$ and $Y$. If $Y$ is a pointed $CW$-complex, then
$|K_\PP|\ast Y$ has the homotopy type of $\Sigma |K_\PP|\wedge Y$.

  {\bf Theorem (Bahri, Bendersky, Cohen and Gitler \cite{BB,BBCG08} 2.12)} {\it Let $K_\PP$ be an abstract simplicial
complex and $\overline{K_\PP}$
its associated poset. Let $(\underline{X}, \underline{A})=\{(X_i, A_i, x_i)|i\in [m]\}$ denote m choices of
connected, pointed pairs of $CW$-complexes, with the inclusion $A_i\hookrightarrow X_i$ homotopic to constant for
all $i$. Then there is a homotopy equivalence
\[
\wh{\ZZ}_{K_\PP}(\underline{X}, \underline{A})\longrightarrow \bigvee_{\omega\in K_\PP}
|\Delta(\overline{K_\PP}_{<\omega})|\ast\wh{D}(\omega).
\]
}

    From (5.2), we see that if $A_i\hookrightarrow X_i$ are null-homotopic for all $i$, then
\begin{align*}
\Sigma(\ZZ_{K_\PP}(\underline{X}, \underline{A}) \simeq &
  \Sigma\left(\bigvee_{I\subset [m]} \wh{\ZZ}_{K_\PP \cap I}(\underline{X_I}, \underline{A_I})\right)\\
   \simeq &  \Sigma\left(\bigvee_{I\subset [m]}\left(\bigvee_{\omega\in K_\PP\cap I}
  |\Delta(\overline{K_\PP\cap I}_{<\omega})|\ast\wh{D}_I(\omega)\right)\right).\notag
\end{align*}
Fix an $\omega\in K_\PP$, it is easy to see that $\omega\in K_\PP\cap I$ if and only if
$\omega\subset I$. Thus $\Sigma(\ZZ_{K_\PP}(\underline{X}, \underline{A})$ has the homotopy type of
\begin{align}\notag
\left(\bigvee_{I\subset [m]}\left(\bigvee_{\omega\in K_\PP\cap I}
  |\Delta(\overline{K_\PP\cap I}_{<\omega})|\ast\wh{D}_I(\omega)\right)\right) = &
\left(\bigvee_{\omega\in K_\PP}\left(\bigvee_{\omega \subset I\subset [m]}
  |\Delta(\overline{K_\PP\cap I}_{<\omega})|\ast\wh{D}_I(\omega)\right)\right)\\
 \simeq & \left(\bigvee_{\omega\in K_\PP}\left(\bigvee_{\omega \subset I\subset [m]}
  \Sigma|\Delta(\overline{K_\PP\cap I}_{<\omega})|\wedge\wh{D}_I(\omega)\right)\right)\tag{\nr}
\end{align}

  {\bf Remark \nr} If $\omega$ is a maximum face of $K_\PP\cap I$
in the sense that $\omega\in K_\PP\cap I$ but it is not a proper subset of any other simplices,
then $\overline{K_\PP\cap I}_{<\omega}=\{\tau\in \overline{K_\PP\cap I}|\omega\varsubsetneq\tau\}=\varnothing$.
The simplicial complex $\Delta(\overline{K_\PP\cap I}_{<\omega})=\{\phi\}$ is the empty simplicial complex.
Here we use the agreement
\begin{align*}
|\{\phi\}|\ast \wh{D}_I(\omega)= & \wh{D}_I(\omega) & \mbox{and} & &
\Sigma|\{\phi\}|\wedge \wh{D}_I(\omega) = & \wh{D}_I(\omega).
\end{align*}
Combine with the agreement in Theorem 2.8, $\w{H}^{-1}(|\{\phi\}|, \kk)=\kk$, we have
\begin{align*}
\w{H}^0(\Sigma|\{\phi\}|, \kk)= & \kk & \mbox{and} & & \w{H}^*(\Sigma|\{\phi\}|\wedge\wh{D}_I(\omega), \kk)= &
\w{H}^*(\wh{D}_I(\omega), \kk)
\end{align*}

   Consider the reduced cohomology of $\ZZ_{K_\PP}(\underline{X}, \underline{A})$, we see from (5.3) that:
\begin{align}\tag{\nr}
\w{H}^*(\ZZ_{K_\PP}(\underline{X}, \underline{A}), \kk)= &
  \bigoplus_{\omega\in K_\PP}\left(\bigoplus_{\omega\subset I\subset[m]}
  \w{H}^*\left(\Sigma|\Delta(\overline{K_\PP\cap I}_{<\omega})|\wedge\wh{D}_I(\omega), \kk\right)\right).
\end{align}

   Furthermore suppose that there is no
$Tor$ problem in the K\"{u}nneth formulae for the cohomology of
$\Sigma|\Delta(\overline{K_\PP\cap I}_{<\omega})|\wedge\wh{D}_I(\omega)$
(For example, take $\kk$ to be a field or suppose that $\w{H}^*(X_i, \kk)$ and $\w{H}^*(A_i, \kk)$ are free
$\kk$-modules for all $i$).
Then from
\begin{align*}
\wh{D}_I(\omega)= & Y_{i_1}\wedge Y_{i_2}\wedge\cdots\wedge Y_{i_k}\simeq
 \left(\bigwedge_{i\in\omega}X_i\right)\bigwedge\left(\bigwedge_{j\in I\setminus\omega} A_j\right),\\
H^*(\wh{D}_I(\omega), \kk) = & \left(\bigotimes_{i\in\omega} \w{H}^*(X_i, \kk)\right) \otimes\left(\bigotimes_{j\in I\setminus\omega}
  \w{H}^*(A_j, \kk)\right)
\end{align*}
which is a free $\kk$-module, we have
\begin{align}\tag{\nr}
 & \w{H}^*\left(\Sigma|\Delta(\overline{K_\PP\cap I}_{<\omega})|\wedge\wh{D}_I(\omega), \kk\right)
=  \w{H}^*\left(\Sigma|\Delta(\overline{K_\PP\cap I}_{<\omega})|, \kk\right)\otimes\w{H}^*\left(\wh{D}_I(\omega), \kk\right)\\
= & \w{H}^*(\Sigma|\Delta(\overline{K_\PP\cap I}_{<\omega})|,\kk)\otimes
  \left(\bigotimes_{i\in\omega} \w{H}^*(X_i, \kk)\right) \otimes\left(\bigotimes_{j\in I\setminus\omega}
  \w{H}^*(A_j, \kk)\right).\notag
\end{align}

 {\bf Theorem \nr} {\it Let $K_\PP$ be an abstract simplicial complex corresponding to a simplicial complement $\PP$
and let $(\underline{X}, \underline{A})=\{(X_i, A_i, x_i)|i\in [m]\}$ denote m choices of
connected, pointed pairs of $CW$-complexes, with the inclusion $A_i\hookrightarrow X_i$ homotopic to constant for
all $i$. If $\w{H}^*(X_i, \kk)$ and $\w{H}^*(A_i, \kk)$ are free $\kk$-modules for all $i$,
then the cohomology  of $\ZZ_{K_\PP}(\underline{X}, \underline{A})$ is isomorphic to
\[
H^*(\ZZ_{K_\PP}(\underline{X}, \underline{A}), \kk)=
\bigoplus_{\omega\in K_\PP}\left(\bigoplus_{\tau}
  H_{*,\tau}(\Lambda^{*,\tau}[E_\omega\PP], d) \otimes
  \left(\bigotimes_{i\in\omega}\w{H}^*(X_i, \kk)\right)\otimes
  \left(\bigotimes_{j\in\tau}\w{H}^*(A_j, \kk)\right)\right)
\]
as $\kk$-modules,  where
\[
H_{*,\tau}(\Lambda^{*,\tau}[E_\omega\PP],d)\cong Tor_{*,\tau}^{\kk[\xx]}(\kk(\mbox{star}_{K_\PP}\omega), \kk)
\]
subject to $\tau\subset [m]\setminus \omega$.
}

\begin{proof}
    Given an $\omega\subset I$, from its definition we see that the posets
\begin{align*}
\overline{K_\PP\cap I}_{\leqslant \omega}= & \{\tau\in K_\PP\cap I|\omega\subset \tau\}
  = \{\tau\in K_\PP|\omega\subset\tau\subset I\} \hspace{10mm}\mbox{and}\\
\overline{\mbox{link}_{K_\PP\cap I}\omega} = & \{\tau'\in K_\PP\cap I|\omega\cup\tau'\in K_\PP\cap I, \tau'\cap\omega=\phi\}\\
  = & \{\tau'\in K_\PP|\omega\cup\tau'\in K_\PP, \tau'\cap\omega=\phi,\ \mbox{and}\ \tau'\subset I\}
     =\overline{(\mbox{link}_{K_\PP}\omega)\cap I}_{\leqslant \phi}.
\end{align*}
There is a one to one correspondence between the posets $\Psi: \overline{K_\PP\cap I}_{<\omega}\longrightarrow
\overline{(\mbox{link}_{K_\PP}\omega)\cap I}_{<\phi}$ given by $\Psi(\tau)=\tau\setminus\omega$.
Thus the order complex $\Delta(\overline{K_\PP\cap I}_{<\omega})=((\mbox{link}_{K_\PP}\omega)\cap I)'$
is the barycentric subdivision of
$(\mbox{link}_{K_\PP}\omega)\cap I$ and then
\begin{align}\tag{\nr}
\w{H}^*(\Sigma|\Delta(\overline{K_\PP\cap I}_{<\omega})|, \kk)=
\w{H}^{*-1}(|\Delta(\overline{K_\PP\cap I}_{<\omega})|, \kk)=
  \w{H}^{*-1}((\mbox{link}_{K_\PP}\omega)\cap I, \kk).
\end{align}

    Recall from Theorem 4.7 that $\mbox{link}_{K_\PP}\omega=K_{E_\omega\PP}\cap([m]\setminus\omega)$, we have
\[(\mbox{link}_{K_\PP}\omega)\cap I =
K_{E_\omega\PP}\cap([m]\setminus\omega)\cap I=K_{E_\omega\PP}\cap (I\setminus\omega).
\]
Then from Theorem 2.8 and (5.8), we have
\begin{align*}
H_{q,I\setminus\omega}(\Lambda^{*,I\setminus\omega}[E_\omega\PP], d)
 =  & \w{H}^{|I\setminus\omega|-q-1}((\mbox{link}_{K_\PP}\omega)\cap I, \kk)
 =  \w{H}^{|I\setminus\omega|-q}(\Sigma|\Delta(\overline{K_\PP\cap I}_{<\omega})|, \kk)
\end{align*}

    From the definition of $E_\omega\PP=\{\sigma_1\setminus\omega, \sigma_2\setminus\omega, \cdots,
\sigma_s\setminus\omega\}$, we see that total subset $S_\uu$ of any generator $\uu\in\Lambda^{*,*}[E_\omega\PP]$
is contained in $[m]\setminus\omega$. Thus the homology of the simplicial complement
$H_{q,\tau}(\Lambda^{*,\tau}[E_\omega\PP],d)$ is concentrated in $\tau\subset[m]\setminus\omega$.
Denote the non-empty set $I$ by $\tau\cup\omega$ for any $\omega\subset I\subset[m]$,
we have
\[
H_{q,\tau}(\lambda^{*,\tau}[E_\omega\PP], d)=\w{H}^{|\tau|-q}(\Sigma|\Delta
 (\overline{K_\PP\cap (\tau\cup\omega)}_{<\omega})|, \kk).
\]

  Apply this formula to (5.6) and (5.5) we get
\begin{align*}
 & \w{H}^*(\ZZ_{K_\PP}(\underline{X}, \underline{A}), \kk)=
  \bigoplus_{\omega\in K_\PP}\left(\bigoplus_{\tau\subset[m]\setminus\omega}
  \w{H}^*\left(\Sigma|\Delta(\overline{K_\PP\cap (\omega\cup\tau)}_{<\omega})|\wedge\wh{D}_{\omega\cup\tau}(\omega), \kk\right)\right) \\
 = &  \bigoplus_{\omega\in K_\PP}\left(\bigoplus_{\tau\subset[m]\setminus\omega}
 H_{q, \tau}(\Lambda^{*, \tau}[E_\omega\PP], d)\otimes
  \left(\bigotimes_{i\in\omega} \w{H}^*(X_i, \kk)\right) \otimes\left(\bigotimes_{j\in \tau}
  \w{H}^*(A_j, \kk)\right)\right)
\end{align*}
where $I=\omega\cup\tau\not=\phi$.

    The theorem follows from $H^0(\ZZ_{K_\PP}(\underline{X}, \underline{A}), \kk)=\kk$ and the agreement by taking
$I=\omega=\tau=\phi$ to be the empty set
\[
H_{0,\phi}(\Lambda^{*,\phi}[E_\phi\PP], d)
  \otimes\left(\bigotimes_{i\in\phi}\w{H}^*(X_i, \kk)\right)
  \otimes\left(\bigotimes_{j\in\phi}\w{H}^*(A_j, \kk)\right)=\kk.
\]
\end{proof}

   {\bf Remark:} If the abstract simplicial complex $K_\PP$ is given by an abstract simplicial complement $\PP$,
we might not know if a subset $\omega$ is a simplex of $K_\PP$ or not. This is not a problem, because
while $\omega\not\in K_\PP$, the empty set $\phi$ is an element of $E_\omega\PP$ and
$H_{*,*}(\Lambda^{*,*}[E_\omega\PP], d)\equiv 0$. Thus Theorem 5.7 could be written as
\begin{align}\notag
  & H^*(\ZZ_{K_\PP}(\underline{X}, \underline{A}), \kk)\\
= & \bigoplus_{\omega\in[m]}\left(\bigoplus_{\tau\in[m]}
  H_{*,\tau}(\Lambda^{*,\tau}[E_\omega\PP], d) \otimes
  \left(\bigotimes_{i\in\omega}\w{H}^*(X_i, \kk)\right)\otimes
  \left(\bigotimes_{j\in\tau}\w{H}^*(A_j, \kk)\right)\right)\tag{\nr}
\end{align}

  {\bf Corollary \nr} {\it If all the $A_i$ are contractible and $\w{H}^*(X_i, \kk)$ are free $\kk$-modules.
Then
\[
H^*(\ZZ_{K_\PP}(\underline{X}, \underline{A}), \kk)=
 \bigoplus_{\omega\in K_\PP}\left(\bigotimes_{i\in \omega} \w{H}^*(X_i, \kk)\right)
\]
}

\begin{proof}
If all the $A_i$ are contractible then $\displaystyle{\bigotimes_{j\in\tau}\w{H}^*(A_j, \kk)=0}$ for any non-empty set $\tau$.
The corollary follows from $H_{0,\phi}(\Lambda^{*,\phi}[E_\omega\PP],d)=\kk$ if $\omega\in K_\PP$.
\end{proof}

 {\bf Corollary \nr} {\it If all the $X_i$ are contractible and $\w{H}^*(A_i, \kk)$ are free
$\kk$-modules. Then
\begin{align*}
H^*(\ZZ_{K_\PP}(\underline{X}, \underline{A}), \kk)= &
 \bigoplus_{\tau}\left(H_{*,\tau}(\Lambda^{*,\tau}[\PP], d)
 \otimes\left(\bigotimes_{j\in\tau}\w{H}^*(A_j, \kk)\right)\right).
\end{align*}
Furthermore take $X_i=D^2$ and $A_i=S^1$ for all $i$,
\begin{align*}
H^{r}(\ZZ_{K_\PP}, \kk)=H^{2|\tau|-q}(\ZZ_{K_\PP}(\underline{D^2}, \underline{S^1}), \kk)= &
 \bigoplus_{2|\tau|-q=r}H_{q,\tau}(\Lambda^{*,\tau}[\PP], d),
\end{align*}
where $H_{q,\tau}(\Lambda^{*,\tau}[\PP], d)=Tor_{q,\tau}^{\kk[\xx]}(\kk(K_\PP), \kk)$.
}

\begin{proof}
    If all the $X_i$ are contractible, then $\displaystyle{\bigotimes_{i\in\omega}\w{H}^*(X_i, \kk)=0}$ for any
non-empty set $\omega$. The corollary follows from
\[
H^*(\ZZ_{K_\PP}(\underline{X}, \underline{A}), \kk)=
 \bigoplus_{\tau}H_{*,\tau}(\Lambda^{*,\tau}[E_\phi\PP], d)
 \otimes\left(\bigotimes_{i\in\phi}\w{H}^*(X_i,\kk)\right)\otimes
  \left(\bigotimes_{j\in\tau}\w{H}^*(A_j, \kk)\right)
\]
and the agreement $\displaystyle{\bigotimes_{i\in\phi}\w{H}^*(X_i,\kk)=\kk}$.
\end{proof}

  {\bf Proposition \nr} {\it Let $X_i=S^2$ and $A_i=S^1$ all $i$. Then
\[
H^{r}(\ZZ_{K_\PP}(\underline{S^2}, \underline{S^1}), \kk) =
  \bigoplus_{\begin{array}{c}\omega\in K_\PP\\ 2|\omega|+2|\tau|-q=r\end{array}}
  \left(\bigoplus_{\tau}H_{q,\tau}(\Lambda^{*,\tau}[E_\omega\PP], d)\right).
\]
}

\begin{proof}
  Noticed that $\wh{D}_{\omega\cup\tau}(\omega)$ is the $2|\omega|+|\tau|$ sphere if $\omega\cap\tau=\phi$,
$$\wh{D}_{\omega\cup\tau}(\omega)
 =\left(\bigwedge_{i\in\omega}S^2\right)\wedge\left(\bigwedge_{j\in\tau}S^1\right)
 =S^{2|\omega|+|\tau|},
$$
we see from (5.5) that
$$H^*(\ZZ_{K_\PP}(\underline{X}, \underline{A}), \kk)=
  \bigoplus_{\omega\in K_\PP}\left(\bigoplus_{\tau\subset[m]\setminus\omega}
  \w{H}^*\left(\Sigma^{2|\omega|+|\tau|+1}|\Delta(\overline{K_\PP\cap (\omega\cup\tau)}_{<\omega}|, \kk\right)\right)
$$
including $\omega=\tau=\phi$. The result follows from
$$\w{H}^{|\tau|-q}\left(\Sigma|\Delta(\overline{K_\PP\cap (\omega\cup\tau}_{<\omega})|, \kk\right)
 =H_{q,\tau}(\Lambda^{*,\tau}[E_\omega\PP], d).$$
\end{proof}

  We finish this paper by giving an example.

{\bf Example \nr} Let $m=6$ and the simplicial complement $\PP=\{\sigma_1=\{1,2\}, \sigma_2=\{3,4\}, \sigma_3=\{5,6\}\}$.
The corresponding simplicial complex is a triangulation of the sphere $S^2$ (see Finger 2)

\begin{center}
\setlength{\unitlength}{0.5mm}
\begin{picture}(100,60)(0,0)
\Thicklines
\drawline(0,0)(100,0)
\drawline(0,0)(0,50)
\drawline(0,0)(25,25)
\drawline(0,50)(25,25)
\drawline(0,50)(100,50)
\drawline(25,25)(125,25)
\drawline(100,0)(125,25)
\drawline(100,50)(125,25)
\thinlines
\drawline(0,0)(125,25)
\drawline(25,25)(100,50)
\dottedline{2}(0,50)(100,0)
\Thicklines
\dottedline{2}(100,0)(100,50)
\put(-6,-1){1}
\put(-6,50){3}
\put(15,22){5}
\put(103,50){2}
\put(127,22){4}
\put(105,-1){6}
\put(35,-20){Finger 2}
\end{picture}
\end{center}
\vspace{15mm}

    The cohomology ring of $\ZZ_K=\ZZ_{K_\PP}(\underline{D^2}, \underline{S^1})$ could be easily got from
the homology of the simplicial complement  $\PP=\{\sigma_1=\{1,2\}, \sigma_2=\{3,4\},$ $\sigma_3=\{5,6\}\}$,
\[
H_{*,*}(\Lambda^{*,\tau}[\PP], d)\cong\Lambda^{*,*}[\{\sigma_1, \sigma_2, \sigma_3\}]
\]
as an exterior algebra over $\kk$. Thus the cohomology ring $H^*(\ZZ_K, \kk)$ is the exterior algebra
on $\PP$ with non-trivial products.
The Poincar\`{e} series of $H^*(\ZZ_K, \kk)$ is  $1+3x^3+3x^6+x^9$ and the
total Betti number of $\ZZ_K$ is $8$.

   Consider the cohomology of $\ZZ_{K_\PP}(\underline{S^2}, \underline{S^1})$. We start from
computing the homology of the simplicial complement $E_\omega\PP$ with $\omega\in K$:
\begin{enumerate}
\item Take $\omega=\phi$,  the homology
of the simplicial complement $E_\phi\PP$ is
$\Lambda^{*,*}[\{\sigma_1, \sigma_2, \sigma_3\}].$
Thus the submodule
\[
\bigoplus_{\tau}H_{*,\tau}(\Lambda^{*,\tau}[E_\phi\PP], d)=\Lambda^{*,*}[\{\sigma_1,\sigma_2,\sigma_3\}]
\]
is a free $\kk$-module with Poincar\`{e} series $1+3x^3+3x^6+x^9$.
\item Take $\omega=\{1\}$ (similarly, take any of its 6 $0$-simplex $\omega=\{i\}$, $i\in[6]$),
the homology
of the simplicial complement $E_{\{i\}}\PP$ is $\Lambda^{*,*}[E_{\{i\}}\PP]$.
Thus the submodules
\[
\displaystyle{\bigoplus_{\tau}H_{*,\tau}(\Lambda^{*,\tau}[E_{\{i\}}\PP], d)}
\]
is the free $\kk$-module with Poincar\`{e} series $x^2(1+x+2x^3+2x^4+x^6+x^7)$.
\item Take $\omega=\{1,3\}$ and similarly any of its 12 $1$-simplex $\omega=\{i_1, i_2\}$, the homology
of the corresponding simplicial complement is $\Lambda^{*,*}[E_{\{i_1,i_2\}}\PP]$. The submodule
\[
\displaystyle{\bigoplus_{\tau}H_{*,\tau}(\Lambda^{*,\tau}[E_{\{i_1, i_2\}}\PP], d)}
\]
is the free $\kk$-module with Poincar\`{e} series $x^4(1+2x+x^2+x^3+2x^4+x^5)$.
\item Take $\omega=\{1,3,5\}$ and similarly any of its 8 $2$-simplex $\omega=\{i_1, i_2,i_3\}$, the homology
of the corresponding simplicial complement is $\Lambda^{*,*}[E_{\{i_1,i_2,i_3\}}\PP]$. The submodule
\[
\displaystyle{\bigoplus_{\tau}H_{*,\tau}(\Lambda^{*,\tau}[E_{\{i_1, i_2\,i_3\}}\PP], d)}
\]
is the free $\kk$-module with Poincar\`{e} series $x^6(1+3x+3x^2+x^3)$.
\end{enumerate}

   Thus
\[
H^*(\ZZ_{K_\PP}(\underline{S^2}, \underline{S^1}), \kk)=
 \bigoplus_{\omega\in K_\PP}\left(\bigoplus_{\tau}H_{q,\tau}(\Lambda^{*,\tau}[E_\omega\PP], d)\right)
\]
is a free $\kk$-module with Poincar\`{e} series
\begin{align*}
  & (1+3x^3+3x^6+x^9)\\
+ & 6x^2(1+x+2x^3+2x^4+x^6+x^7) \\
+ & 12x^4(1+2x+x^2+x^3+2x^4+x^5)\\
+ & 8x^6(1+3x+3x^2+x^3)\\
= & 1+6x^2+9x^3+12x^4+36x^5+35x^6+36x^7+54x^8+27x^9
\end{align*}
The total Betti number of $\ZZ_{K_\PP}(\underline{S^2}, \underline{S^1})$ is $216$.

\end{document}